\documentclass[11pt,reqno]{article}
\usepackage[utf8]{inputenc}
\usepackage{bm, amsmath, mathrsfs, amssymb, amsthm, amscd, amsfonts, graphicx, color,epsfig,latexsym,eucal,layout,fancyhdr, colonequals, enumerate, hyperref, fontenc, hyperref, graphics, a4wide, verbatim, psfrag}

\newcounter{EQNR}
\renewcommand{\contentsname}{Contents\\
{\footnotesize\normalfont(A table of contents should normally not be included)}}

\newtheorem{theorem}{Theorem}
\newtheorem{corollary}[theorem]{Corollary}

\newtheorem{example}[theorem]{Example}
\newtheorem{lemma}[theorem]{Lemma}
\newtheorem{proposition}[theorem]{Proposition}

\newtheorem{remark}[theorem]{Remark}

\begin{document}

\title{On an approach for evaluating certain trigonometric character sums using the discrete time heat kernel}
\author{Carlos A. Cadavid,\and Paulina Hoyos, \and Jay Jorgenson,\footnote{The third-named author acknowledges grant support
from PSC-CUNY Award 65400-00-53, which was jointly funded
by the Professional Staff Congress and The City University of New York.} \and
Lejla Smajlovi\'{c}, \and Juan D. V\'elez}
\maketitle

\begin{abstract}\noindent
In this article we develop a general method by which one can explicitly
evaluate certain sums of $n$-th powers of products of $d\geq 1$ elementary trigonometric
functions evaluated at $\mathbf{m}=(m_1,\ldots,m_d)$-th roots of unity.  Our approach is to first identify
the individual terms in the expression under consideration as eigenvalues
of a discrete Laplace operator associated to a graph whose vertices form a $d$-dimensional discrete torus $G_{\mathbf{m}}$
which depends on $\mathbf{m}$.  The sums in question are then related to
the $n$-th step of a Markov chain on $G_{\mathbf{m}}$.  The Markov chain admits
the interpretation as a particular random walk, also viewed as a discrete time and discrete space heat
diffusion, so then the sum in question is related to special values of the associated
heat kernel.  Our evaluation follows by deriving a combinatorial expression for
the heat kernel, which is obtained by periodizing the
heat kernel on the infinite lattice $\mathbb{Z}^{d}$ which covers $G_{\mathbf{m}}$.
\end{abstract}

\section{Introduction}

\subsection{A motivating example}

There are some intriguing trigonometric identities which would catch the attention
of almost any mathematician, or mathematics student for that matter.  For example,
for any positive integers $n$ and $m$, we have that
\begin{equation} \label{eq: example1}
\sum\limits_{j=0}^{m-1}\cos ^{n}\left(\frac{2\pi j}{m}\right)=2^{-n}m\sum\limits_{k=-%
\left\lfloor n/m\right\rfloor ,\text{ \ }km+n\text{ even}}^{\left\lfloor
n/m\right\rfloor } \binom{n}{(km+n)/2},
\end{equation}
where $\lfloor x\rfloor$ stands for the largest integer less than or equal to $x$.
%and $\binom{n}{h}$ is a binomial coefficient.
Of course, for small $m$, the values
of $\cos(2\pi j/m)$ can be explicitly computed.  For instance, if $m=5$ then $$\{\cos(2\pi j/5)\}_{j=1}^{4} = \{(\pm 1 \pm \sqrt{5})/4\}.$$
With this, the verification of \eqref{eq: example1} for small $n$ becomes (or at least should
be) an enjoyable exercise for students in high school mathematics courses.  The analogous
expression in the case of $m=7$ is even more pleasant, though challenging.

As is often true in mathematics, there is much more to \eqref{eq: example1} than the
formula itself.  For example, if one multiplies both sides of \eqref{eq: example1} by $2^{n}$,
the left-hand-side is
\begin{equation}\label{eq:example1_part2}
\sum\limits_{j=0}^{m-1}\left(e^{2\pi i j/m} + e^{-2\pi ij/m}\right)^{n}.
\end{equation}
Each summand in \eqref{eq:example1_part2} is an algebraic integer in the cyclotomic field $\mathbb{Q}(\zeta_{m})$ where
$\zeta_{m}$ is a primitive $m$-th root of unity.  It is elementary to prove that the sum \eqref{eq:example1_part2}
is invariant under the action of the Galois group $\text{\rm Gal}(\mathbb{Q}(\zeta_{m})/\mathbb{Q})$;
hence, \eqref{eq:example1_part2} is an integer in $\mathbb{Q}$.  Therefore, \eqref{eq: example1}
can be viewed as giving
an explicit evaluation of the rational integer \eqref{eq:example1_part2}.  Indeed,  one can obtain
\eqref{eq: example1} from \eqref{eq:example1_part2} by applying the binomial theorem.

\subsection{The general question}

The proof of \eqref{eq: example1} described above is specific
to the series in question and
does not admit a complete understanding of the mathematical structure behind similar questions.
With this assertion in mind, the purpose of this article (and the follow-up article \cite{CHJSV22a})
is to address the following general problem.

\it Let $d\geq 1$ be an
integer, and let $P(x_{1},\cdots, x_{d};y_{1},\cdots, y_{d})$ be any polynomial in $2d$ variables with complex coefficients.
Let $\{m_{j}\}_{j=1}^{d}$ and $\{a_{j}\}_{j=1}^{d}$ be sets of positive integers, and
$\{\beta_{j}\}_{j=1}^{d}$ a set of real numbers.  Can one determine an efficient
and effective algorithm which evaluates
\begin{multline}\label{eq:General_sum}
\sum\limits_{k_{1}=0}^{m_{1}-1} \cdots \sum\limits_{k_{d}=0}^{m_{d}-1}
P\left(\cos\left(\frac{2\pi k_{1}a_{1}}{m_{1}}+ \beta_{1}\right),\cdots , \cos\left(\frac{2\pi k_{d}a_{d}}{m_{d}}+ \beta_{d}\right); \right.\\ \left.\sin\left(\frac{2\pi k_{1}a_{1}}{m_{1}}+ \beta_{1}\right),\cdots , \sin\left(\frac{2\pi k_{d}a_{d}}{m_{d}}+ \beta_{d}\right) \right)
\end{multline}
as a finite sum involving the degree and coefficients of $P$, binomial coefficients
involving the given data, and exponentials
which are linear in $\beta_{j}$, $j=1,\ldots,d$? \rm

Specific instances of \eqref{eq:General_sum}, as well as similar series which include additive or multiplicative characters, or trigonometric realization of important sums, such as Dedekind or Hardy sums
have attracted considerable interest, in part because such series appear in different mathematical and physical settings; see, for example, \cite{Ve88},
\cite{Do89}, \cite{Do92}, \cite{CS12} and \cite{Ho18} for sums within physical settings or \cite{MS20}, \cite{Ro19}, \cite{dFGK18} and \cite{So18} for a
other mathematical research on related sums.  As a result, one can see the potential significance in obtaining
closed-form evaluations of such sums.  However, as stated in \cite{BY02},
many of those sums seem intractable or simply do not have known evaluations, however, it is possible to establish reciprocal relations for such sums;
see also  \cite{BC13}, \cite{Ch18} and \cite{MS20}.

Various approaches have been used to evaluate certain finite trigonometric sums.  Here is a partial list
of articles we found particularly interesting.

In \cite{BY02}, the authors used contour integration, and the article includes an
interesting discussion of the history of this method.  In \cite{CM99}, certain series were evaluated
using a generating series and a partial fraction decomposition.  In \cite{AH18}, the authors proved numerous relations by appealing to results in the theory of special functions.
In \cite{Ha08}, some trigonometric sums were computed by using a discrete form of the sampling theorem associated with certain second-order discrete eigenvalue problems. In \cite{AZ22} the authors show that many trigonometric identities have been ``re-discovered'' many times in the vast literature on the topic and describe an interesting ``automated approach'' for proving some types of trigonometric identities.

Each of the above mentioned
approaches has proved interesting formulas in specific instances or with certain types of series.  Still, it
would be of mathematical value to develop any other means by which
one can compute \eqref{eq:General_sum} in closed form.

\subsection{Our approach}

The aim of this note is to present a unified approach to the numerical evaluation of a wide class of finite trigonometric sums
involving multiple sums of products of powers of sine and cosine functions evaluated at rational multiples of $\pi$, possibly
twisted by a character or shifted by an arbitrary additive constant.  In effect, we consider the setting of
\eqref{eq:General_sum} in the case $P$ is a monomial.

More precisely, we express the discrete time heat kernel on a $d$-dimensional discrete torus $G_{\mathbf{m}}$, twisted by a character of $\mathbb{Z}^d$ in two different ways: First, through a combinatorial approach, stemming from the
realization of $G_{\mathbf{m}}$ as a quotient of $\mathbb{Z}^{d}$ by a sublattice; and second, through a spectral theoretic approach,
which comes from the expansion of the heat kernel in terms of eigenfunctions and eigenvalues of the graph Laplacian.
The eigenvalues of the discrete Laplacian are, in effect, the trigonometric terms we seek to evaluate, and the identities
follow from the general result that the discrete time heat kernel is unique, hence any two evaluations are equal.
It is important to note that we consider the setting where the edges of $G_{\mathbf{m}}$ are not just ``nearest neighbor'',
which allows us considerable flexibility in our analysis.

Our approach is similar in spirit to the method undertaken by Ejsmont and Lehner in \cite{EL21} who also obtain an evaluation
of a trigonometric sum by employing spectral theory. For instance,
Corollary  \ref{cor: main trace f-la} below is proved by looking at the trace of a discrete time heat kernel. However,
our approach is more general since we are looking at the complete spectral expansions in $d$ dimensions instead of simply
the trace of an operator alone.
In doing so, we are able to consider many trigonometric sums with a character twist.

Additionally, our heat kernel approach
provides a physical interpretation of many finite trigonometric sums twisted by a character as a spectral expansion
on a specially designed weighted Cayley graph.
For example, for positive integers $m_1,m_2, a,b, k$ and real
numbers $\alpha_1,\alpha_2$ the trigonometric sum
\begin{equation}\label{eq:sample}
\sum_{j=0}^{2m_1-1}\sum_{\ell=0}^{m_2-1}(-1)^j\cos^k\left(\frac{\pi ja}{m_1}+\alpha_1\right)\sin^k\left(\frac{2\pi \ell b}{m_2}+\alpha_2\right)
\end{equation}
equals $2m_1m_2\exp(-\pi i \alpha_1 m_1/a)$ times a specific value of the heat kernel associated to the weighted Cayley graph
$\mathcal{C}(\mathbb{Z}^2/\mathbf{m}\mathbb{Z}^2,S,\pi_S)$ where $\mathbf{m}=(2m_1,m_2)$, $S=\{(\pm a, \pm b)\}$ is a set of four elements
and $\pi_S(s)=1/4$ for all $s\in S$.  Furthermore, the discrete Laplacian is twisted by a character $\chi_{\bm{\beta}}$ with $\bm{\beta}=\left(\frac{\alpha_1 m_1}{\pi a}, \frac{m_2(\alpha_2-\pi/2)}{2\pi b} \right)$.  With all this, the heat kernel in question is evaluated at the space variables
$\mathbf{x}=(m_1,0)$, $\mathbf{y}=(0,0)$ and
at discrete time $k$, from which we can evaluate \eqref{eq:sample}.
For additional details, see Section \ref{sec: discrete time HK on discrete tori} which gives a detailed explanation of the notation.
The combinatorial evaluation of the discrete time heat kernel, as given in Proposition \ref{ heat kernel geometric sum},
combined with Proposition \ref{prop: product of HK} yields an explicit and closed evaluation of this sum; namely, it is proved in Example \ref{ex: final d=2}
that
\begin{equation} \label{eq:example 20}
\eqref{eq:sample}= \frac{2m_1m_2}{4^k} e^{ik(\alpha_1+\alpha_2)}\cdot \sum\limits_{_{\substack{ d_1,d_2 \in\{0,\ldots, k\} \\ a(2d_1-k)/m_1-1/2 \in \mathbb{Z}, \text{  } b(2d_2-k)/m_2 \in \mathbb{Z} }}} \binom{k}{d_1}  \binom{k}{d_2}  i^{2d_2-k} e^{-2i(\alpha_1d_1+\alpha_2d_2)}.
\end{equation}
In other words, by first looking at a series in question, as in \eqref{eq:sample}, one can consider a type of
 ``reverse engineering" to construct the graph from which our approach will evaluate the sum.  For example,
in Section \ref{further_sums} below we describe the setting from which one will obtain an explicit evaluation of the sum
\begin{equation}\label{eq:sum_powers}
\sum\limits_{k_{1}=0}^{m_{1}-1}\cdots
\sum\limits_{k_{d}=0}^{m_{d}-1}\left(\cos^{h_{1}}\left(\frac{2\pi k_{1}a_{1}}{m_{1}}+ \beta_{1}\right)+\cdots +
\cos^{h_{d}}\left(\frac{2\pi k_{d}a_{d}}{m_{d}}+ \beta_{2}\right)\right)^{n}.
\end{equation}
For the sake of space, we do not write the evaluation;
however, we are able to give a reasonably brief yet thorough description as to how one can apply
our general result in order to obtain a closed form evaluation of \eqref{eq:sum_powers} in terms of
binomial coefficients and exponentials.

\subsection{Outline of the paper and overview of results}

This article is organized as follows.  In Section 2 we present background material from the literature and
establish the notation and language needed throughout the article.  The main geometric object of study is a Cayley graph $X=\mathcal{C}(G,S,\pi_S) $ associated to a finite abelian group $G$ with a reasonably general set of edges generated by $S\subseteq G$ with corresponding
edge weights $\pi_{S}$.  The only requirement
on $S$ is that it is symmetric, so then two vertices are connected by an edge if and only if their difference belongs to $S$.
An edge $s$ connecting $x,y \in G$ where $x-y=s\in S$ has associated weight $\pi_S(s)$.
The edge weights $\pi_{S}$ are assumed to be positive.  For every vertex, the sum of all weights of the joining edges is one.
Therefore, the weights can be viewed as defining a
probability distribution $\pi_S$ on $S$.  For this article, we take $G$
to be a discrete torus, and we write $G_{\mathbf{m}}=\mathbb{Z}^d/\mathbf{m}\mathbb{Z}^d$.
Also, we can identify $G_{\mathbf{m}}$ as the cosets of the action of the group of translations by a vector
$\mathbf{m}=(m_1,\ldots, m_d) \in \mathbb{Z}^d$ on the lattice $\mathbb{Z}^d$. Of course, $m_{j}>0$ for all $j$.  For a Cayley graph $X$, one has a naturally
defined Laplacian matrix, and the heat kernel is another operator defined from the Laplacian matrix.
%Again, the necessary notation and language so that this article is self-contained is given in Section 2.

In Section 3, we consider fixed but arbitrary $\mathbf{m}$, $S$ and $\pi_S$.  We then
derive one of the two evaluations of the discrete time heat kernel on $X=\mathcal{C}(G_{\mathbf{m}},S,\pi_S) $.
We begin by evaluating the discrete time heat kernel on the Cayley graph $Y=\mathcal{C}(\mathbb{Z}^{d},S,\pi_S)$ as a single combinatorial
coefficient, from which we prove that the discrete time heat kernel on $X$ can
be expressed in terms of a sum of combinatorial coefficients; see Theorem \ref{thm:hk_formula_quotient}.
The heat kernel evaluated in Lemma \ref{lem:HK on Zd} and Theorem \ref{thm:hk_formula_quotient} may be viewed as a solution to a
discrete time dynamic diffusion equation on the $d$-dimensional lattice $\mathbb{Z}^d$ and the discrete torus $G_{\mathbf{m}}$, respectively, where the coefficients of the equation are given by the probability distribution $\pi_S$.

The results of Section 3 can also be viewed as a generalization of discrete time results of \cite{SS14} and \cite{SS15}
 to the case when the space variable belongs to the $d$-dimensional lattice or a $d$-dimensional discrete torus. Also, the papers
 \cite{SS14} and \cite{SS15} consider the diffusion process related to the weighted Laplacian on the graph $Y=\mathcal{C}(\mathbb{Z}^{d},S,w_S)$,
 where $w_S$ is more general, not necessarily non-negative weight function on the symmetric subset $S$ of $\mathbb{Z}$.
 Though we are dealing with probability distribution only, we assert that our methods can be adopted to treat the general setting which
 involve non-negative weights.

In Section 4 we prove the main theorem of the article, which is Theorem \ref{thm: main}, by obtaining a second expression
for the discrete time heat kernel.  As stated above, the second formula for the discrete time heat kernel
follows form the spectral theorem of the associated Laplacian matrix.  In Section 5 we develop several examples of Theorem \ref{thm: main}
which include reproving some results in the literature by utilizing special instances of our main theorem.  Also,
we obtain a new proof of the evaluation of classical Gauss sums.

Of the examples given, we find Corollary \ref{cor:NewNumerics} particularly intriguing.
Our result yields a new combinatorial expression for the trigonometric sum, which is evaluated in the main theorem of \cite{dFK13}. Specifically,
the evaluation of $S(100,13)$ in \cite{dFK13} follows from their algorithm, whereas our method yields a specific formula.  The
numerical expressions, which of course must coincide, are explicitly given and certainly are not a value which anyone would guess.
Going further, in Section 5, Corollary \ref{cor: multipl character sums} we derive two explicit evaluations of sums of powers of
sines and cosines twisted by a primitive (multiplicative) Dirichlet character, related
to results of \cite{BH10} (Theorem 1.2.), \cite{BZ04}, \cite{BBCZ05} and \cite{ZW18}.

In Section 6 we prove a general formula which amounts to showing that the discrete time heat kernel on the
product of Cayley graphs is equal to the corresponding product of discrete time heat kernels.  In some sense, one
could view the result as ``contained in the mathematical folklore''.  However, for the sake of completeness in this article as well as for future studies, the mathematical statement is presented and proved.
As an application of the general formula in Section 6 we prove \eqref{eq:example 20} and deduce some further special formulas, such as the identity
\begin{eqnarray} \label{eq product of d terms}
		&&\sum\limits_{\ell _{1}=0}^{m_{1}-1}\cdots \sum\limits_{\ell
			_{d}=0}^{m_{d}-1}\cos ^{n}\left(\frac{2\pi (\ell _{1}+\beta _{1})}{m_{1}}\right)\cdots
		\cos ^{n}\left(\frac{2\pi (\ell _{d}+\beta _{d})}{m_{d}}\right) \\
		&=&2^{-dn}m_{1}\cdots m_{d}\sum\limits_{k_{1}=-\left\lfloor
			n/m_{1}\right\rfloor }^{k_{1}=+\left\lfloor n/m_{1}\right\rfloor }\cdots
		\sum_{k_{d}=-\left\lfloor n/m_{d}\right\rfloor }^{k_{d}=+\left\lfloor
			n/m_{d}\right\rfloor }e^{-2\pi i(k_{1}\beta _{1}+\ldots k_{d}\beta _{d})}%
		\binom{n}{\frac{n+k_{1}m_{1}}{2}}\cdots \binom{n}{\frac{n+k_{d}m_{r}}{2}}, \nonumber
	\end{eqnarray}
	where the sum on the right-hand side of the above equation is taken only
	over $k_{j}$ for which the sum $n+k_{j}m_{j}$, $j=1,\ldots ,d$ is even.

Finally, in Section 7, we offer a concluding remark in which we describe how to use the results given in detail in this
article in order to circle back and evaluate \eqref{eq:General_sum}.  As the discussion shows, it is vital that
our previous results are proved by allowing for the real character twist.

%The remarks in Section 7 are a short, albeit complete, outline of the approach, which we will present in detail in \cite{CHJSV22a}.

In summary, we devise a complete and general approach by which one can evaluate sums of the form \eqref{eq:General_sum}
in terms of the data associated to the polynomial $P$, such as its degree and coefficients, in terms of
combinatorial coefficients.  Each specific instance of the general result yields an identity, and each identity
is just as enjoyable as \eqref{eq: example1}.  However, formulas such as \eqref{eq: example1} are not isolated, but
rather part of the general mathematical structure coming from the discrete time heat kernel associated to the Cayley graph of finite
abelian groups.

\subsection*{Acknowledgment}

We thank the anonymous referees for many valuable suggestions which helped us improve the clarity of the exposition.

\section{Preliminaries}

\subsection{Weighted Cayley graphs of abelian groups}
Let $G$ be a finite or countably infinite additive abelian group. %If $G$ is countably infinite, we assume $G$ to be countable.
Let $S\subseteq G$ be a finite symmetric subset of $G$. The symmetry condition means that if $s\in S$ then $-s\in S$. We will not
assume that $S$ generates $G$. Let $\alpha \colon S \to \mathbb{R}_{>0}$ be a function such that $\alpha(s) = \alpha(-s)$. The
weighted Cayley graph $X= \mathcal{C}(G, S, \alpha)$ of $G$ with respect to $S$ and $\alpha$ is constructed as follows. The vertices
of $X$ are the elements of $G$, and two vertices $x$ and $y$  are connected with an edge if and only if $x-y \in S$.  The weight
$w(x,y)$ of the edge $(x,y)$ is defined to be $w(x,y) \colonequals \alpha(x-y)$.

One can show that $X$ is a regular graph of degree
\[d= \sum_{s \in S} \alpha (s).\]
If $\alpha$ is a probability distribution on $S$, then the degree of the graph $X$ equals $1$. In this case, we will denote $\alpha$ by $\pi_S$.

A function $f: G\to \mathbb{C}$ is an $L^2$-function if $\sum_{x\in G} |f(x)|^2 <\infty$. The set of $L^2$-functions on $G$
is a Hilbert space $L^2(G, \mathbb{C})$ with respect to the classical scalar product of functions
\[ \langle f_1,f_2 \rangle = \sum_{x\in G} f_1(x)\overline{f_2(x)}.\]

We will denote by $\delta_x$ the standard delta function given by $\delta_x(x) = 1$ and $\delta_x(y) = 0$ for $x \neq y$.

The adjacency operator  $A_X \colon L^2(G, \mathbb{C}) \to L^2(G, \mathbb{C})$ of the graph $X$ is defined as
$$
(A_X f)(x) = \sum_{x-y \in S} \alpha(x-y)f(y).
$$
When $X$ is finite, the adjacency operator written in the standard basis is called the adjacency matrix $A_X$ of the graph $X$. The $(x,y)$-entry of the adjacency matrix is $A_X(x,y) = \alpha(x-y)$. Since $\alpha(x-y)=\alpha(y-x)$, the matrix $A_X$ is symmetric.
Moreover, when $\alpha = \pi_S$, it has the property that the elements of its columns, and rows since $A_X$ is symmetric, sum up to $1$.

Given $x \in G$, let $\chi_x$ denote the character of $G$ corresponding to $x$; see, for example, \cite{CR62} and
Section \ref{sec: discrete time HK on discrete tori} below.
As proved in Corollary 3.2 of \cite{Ba79},
the character $\chi_x$ is an eigenfunction of the adjacency operator $A_X$ of $X$ with corresponding eigenvalue
\begin{equation*}
\eta_{x}=\sum_{s\in S} \alpha(s) \chi_{x}(s).
\end{equation*}

%If the group $G$ is infinite, countable abelian group, for a given symmetric set $S$ that generates $G$ and a positive probability distribution $\pi_S$ on $S$, the Cayley graph $X$ is well-defined as above and it is a regular graph of degree one.

\subsection{Discrete time heat kernel on weighted Cayley graphs}

Let $X$ denote the weighted Cayley graph $\mathcal{C}(G,S,\pi_S)$. The standard, or random walk,
Laplacian $\Delta_X$ is defined to be the operator on $L^2(G, \mathbb{C})$ given by
\[\Delta_X f(x) = f(x)-\sum_{x-y \in S} \pi_S(x-y)f(y).
%= f(x)-\sum_{x\sim y}\pi_S(x-y)f(y).
\]
When $G$ is finite, the random walk Laplacian can be represented by the matrix $I-A_X$.

%The discrete time heat kernel $K_X: X\times \mathbb{Z}_{\geq 0} \to \mathbb{R}$ on $X$ is then defined as a \textcolor{red}{unique (either prove or give reference)} solution to the discrete-time heat equation \[(\Delta_X+\partial_n)K_X(x;n)=0 \] with initial condition $K_X(x;0)=p_0(x)$, $x\in G$, for a given probability distribution $p_0$ on $G$.

%As is common in the literature,
The discrete time heat kernel $K_X: X\times X \times \mathbb{Z}_{\geq 0} \to \mathbb{R}$ on $X$ %on the Cayley graph $X$
is defined as the unique solution of the equation
\begin{equation} \label{eq: discrete heat eq inf graph}
(\partial_{n} + \Delta_X) K_X(x,y;n)= 0, \quad n\geq 0,
\end{equation}
viewed as a function of $x\in G$ for a fixed $y\in G$, and with initial condition $K_X(x,y;0)=\delta_{x}(y)$. It can be shown that this also holds if we interchange the roles of $x$ and $y$. Here $\partial_n$ denotes the discrete time derivative, which is defined as
\[\partial_n K_{X}(x,y;n)=K_{X}(x,y;n+1)- K_{X}(x,y;n).\]

A straightforward computation shows that \eqref{eq: discrete heat eq inf graph} is equivalent to
\[K_X(x,y;n+1)= \sum_{z - x \in S}\pi_S(x-z)K_X(z,y;n) = A_X^{n+1} K_X(x,y;0),\]
as a function of $x$. Then, when $G$ is finite, $ K_X(x,y;n)$ corresponds to the $(x,y)$-entry of the matrix $A_X^n$. More intrinsically, $K_X(x,y;n)$ corresponds to $A_X^n(\delta_y)(x).$ Moreover, one can also show that
$K_X(x,y;n)=K_X(x-y,0;n)$ when $x$ and $y$ are in the same connected component. % and $0_G$ denotes the neutral element of the additive abelian group $G$.

%The heat kernel $K_X(x,y;n)$ can be interpreted as the probability that a particle that starts at $y$ ends at $x$ after $n$ steps. We assume that a particle at vertex $z$ moves to an adjacent vertex $w$ with probability $\pi_S(z-w)$.

The heat kernel $K_X(x,y;n)$ can be interpreted as the probability that the random walk, starting at $y$ arrives at $x$ after $n$ steps. For this random walk we assume that a particle at  $z$ moves to an adjacent vertex $w$ with probability $\pi_S(z-w)$.

\subsection{Twisted discrete time heat kernel on discrete tori}\label{sec: discrete time HK on discrete tori}
For fixed integers $d\geq 1$, $n > 0$, and a vector $\mathbf{m}=(m_1,\ldots,m_d)\in\mathbb{Z}_{>0}^d$, we denote by $G_{\mathbf{m}}$  the group $\mathbb{Z}^d/\mathbf{m}\mathbb{Z}^d = \mathbb{Z}/m_1\mathbb{Z} \times \cdots \times \mathbb{Z}/m_d\mathbb{Z}$, which is a discrete torus, i.e. a product of $d$ discrete circles. Let $S$ be an arbitrary (not necessarily generating) symmetric subset of $G_{\mathbf{m}}$ and $\pi_S$ a probability distribution on $S$. Let  $X = \mathcal{C}(G_{\mathbf{m}},S, \pi_S)$ be  the corresponding Cayley graph.

For an arbitrary vector $\bm{\beta}=(\beta_1,\ldots, \beta_d)$ of real numbers, consider the character $\chi_{\bm{\beta}}:
\mathbb{Z}^d \to \mathbb{C}$ of $\mathbb{Z}^d$ defined as $\chi_{\bm{\beta}}(\mathbf{x})=\exp(2\pi i \bm{\beta}\cdot \mathbf{x})$. Under this definition, it suffices to take $\bm{\beta} \in \mathbb{R}^d/\mathbb{Z}^d$.

We are interested in the space $L^2(\mathbb{Z}^d, \mathbb{C},\chi_{\bm{\beta}})$ of $L^2$-functions $f: \mathbb{Z}^d \to \mathbb{C}$ satisfying the transformation property
\begin{equation}\label{eq:twist by character transf prop}
f(\mathbf{x}+\mathbf{k}\mathbf{m})=\exp(2\pi i \bm{\beta}\cdot \mathbf{k})f(\mathbf{x})
\end{equation}
under the (additive) action of the group $m_1\mathbb{Z} \times \cdots \times m_d\mathbb{Z} =\mathbf{m}\mathbb{Z}^d$.
Here we use the notation  $\mathbf{k}\mathbf{m}= (k_1m_1,\ldots, k_dm_d)$, for $\mathbf{k}=(k_1,\ldots,k_d)$ and $\mathbf{m}=(m_1,\ldots,m_d)$, and $\bm{\beta}\cdot \mathbf{k} = \beta_1k_1 + \ldots + \beta_d k_d$ for the classical dot product.

A function $f$ satisfying \eqref{eq:twist by character transf prop} will be called a
$\chi_{\bm{\beta}}$-twisted $\mathbf{m}$-periodic function on $\mathbb{Z}^d$. When $\bm{\beta}=(0,\ldots,0)
\in \mathbb{Z}^d$, the function $f$ is periodic with period $\mathbf{m}$. When $\bm{\beta}=(1/2,\ldots,1/2)$,
functions $f$ satisfying transformation property \eqref{eq:twist by character transf prop} are sometimes called anti-periodic.

%The probability distribution $\pi_S$ guides a random walk or a diffusion process on $\mathcal{C}(G_{\mathbf{m}},S)$ in the following way: for each vertex $v$ of $X$, the probability of moving to an adjacent vertex $w$ equals $\pi_S(v-w)$ and the probability of moving to a non-adjacent vertex equals zero. In other words, this random walk is a Markov chain with state space $(\mathcal{C}(G_{\mathbf{m}},S),\mathcal{P}(\mathcal{C}(G_{\mathbf{m}},S)))$, arbitrary initial probability distribution $p_0$ and transition probability matrix equal to the adjacency matrix $A$ of the weighted graph   $X=\mathcal{C}(G_{\mathbf{m}},S, \pi_S)$. The random walk (or diffusion) normalized Laplacian is defined as $\Delta=I-A$, where $I$ is the identity matrix.
\vskip .06in

The \it $\chi_{\bm{\beta}}$-twisted discrete time heat kernel \rm
on $X$ is defined to be a function
\begin{equation}\label{eq:beta_twisted_HK}
K_{X,\bm{\beta}}(\mathbf{x},\mathbf{y};n): \mathbb{Z}^d \times \mathbb{Z}^d \times\mathbb{Z}_{\geq 0} \to \mathbb{R},
\end{equation}
and it has the following properties.
For a fixed $\mathbf{y} \in \mathbb{Z}^d$, and viewed as  a function of $\mathbf{x}\in \mathbb{Z}^d$, \eqref{eq:beta_twisted_HK} satisfies
\eqref{eq:twist by character transf prop}, and also as a function of $\mathbf{y}\in \mathbb{Z}^d$ for a fixed $\mathbf{x}\in \mathbb{Z}^d$.
Additionally, when viewed as a function of  $n$, \eqref{eq:beta_twisted_HK}
satisfies the discrete analogue of the heat, or diffusion, equation
$$
(\Delta_X + \partial_n)K_{X,\bm{\beta}}(\mathbf{x},\mathbf{y};n) =0,
$$
with initial condition $K_{X,\bm{\beta}}(\mathbf{x},\mathbf{y};0)=\delta_{\mathbf{y}}(\mathbf{x})$, for $\mathbf{x},\mathbf{y} \in \mathbb{Z}^d/\mathbf{m}\mathbb{Z}^d$.

\section{A combinatorial formula for the twisted  heat kernel}

In this section  we will study \eqref{eq:beta_twisted_HK}, which is
the discrete time heat kernel  on the Cayley graph $X = \mathcal{C}(G_{\mathbf{m}},S, \pi_S)$ twisted
by $\chi_{\bm{\beta}}$. The heat kernel on $X$ can be computed from the heat kernel $K_Y(\mathbf{x}, \mathbf{y};n)$ on the Cayley
graph $Y=\mathcal{C}(\mathbb{Z}^d, S,\pi_S)$ by using the methods which are common in the theory of automorphic forms,
and is analogous to the construction of the continuous time heat kernel on the discrete torus; see \cite{KN06}, \cite{CJK10} and \cite{Do12}.

We start by determining $K_Y(\mathbf{x}, \mathbf{y};n)$. Since $K_Y(\mathbf{x}, \mathbf{y};n)= K_Y(\mathbf{x}-\mathbf{y},  \mathbf{0};n)$,
it suffices to compute $ K_Y(\mathbf{x},  \mathbf{0};n)$. A general formula is presented in the lemma below.

\begin{lemma} \label{lem:HK on Zd} Let $S$ be a symmetric, finite subset of $\mathbb{Z}^d$ and let $\pi_S$ a probability distribution on $S$. For $j\in \{1,\ldots, d\}$ let $S_j:=\max\{s_j:s=(s_1,\ldots,s_d)\in S\}.$
By $\mathbf{t}$ we denote the $d$-tuple $(t_1,\ldots,t_d),$ and for $\mathbf{q}=(q_1,\ldots,q_d)\in \mathbb{Z}^d$ we let $\mathbf{t}^{\mathbf{q}}$ denote $t_1^{q_1}\ldots t_d^{q_d}$.
\begin{itemize}
\item[a.]
If $|x_j|>nS_j,$ for some $j\in \{1,\ldots, d\}$, then
$$ K_Y(\mathbf{x},  \mathbf{0};n)= K_Y((x_1,\ldots,x_d), \mathbf{0};n) =0.$$
\item[b.]
If  $x_j\in\{-nS_j,\ldots,0,\ldots,nS_j\}$, for all  $j\in \{1,\ldots, d\}$, then $ K_Y(\mathbf{x},
\mathbf{0};n)$ equals the coefficient multiplying $\mathbf{t}^{\mathbf{x}}$ in the expansion of the polynomial
\begin{equation}
	\label{sum of monomials}
	\left(\sum_{s\in S} \pi_S(s) \mathbf{t}^{s}\right)^n.
\end{equation}
\end{itemize}
\end{lemma}
\begin{proof}
First, notice that $S_j\geq 0$ because $S$ is assumed to be symmetric;  hence, the conditions given on coordinates $x_j$ of
$\mathbf{x}$ are well-defined.
The proof of the lemma follows immediately by induction, using the observation that, for a fixed $\mathbf{x}\in \mathbb{Z}^d$,
$ K_Y(\mathbf{x},  \mathbf{0};n)$ is the probability that the random walk, starting at $ \mathbf{0}$, ends up at $\mathbf{x}$
after $n$ steps. Therefore, at each step, a particle at a point $\mathbf{y} \in \mathbb{Z}^d$ can move only to points
$\mathbf{y} +s$, for $s\in S$ with probability $\pi_S(s)$.
\end{proof}

In the proposition below, the heat kernel $K_{X}(\mathbf{x},\mathbf{y};n)$ twisted by a character $\chi_{\bm{\beta}}$ is derived from the heat kernel on $Y$ using what could be called the \it method of twisted images; \rm see, for example, \cite{Do12}.

\begin{proposition} \label{ heat kernel geometric sum}
  With notation as above, the discrete time heat kernel on $X$ twisted by the character $\chi_{\bm{\beta}}$ on $\mathbb{Z}^d$ is given by
  \begin{equation} \label{eq: heat k on tori as a sum}
  K_{X,\bm{\beta}}(\mathbf{x}, \mathbf{y};n)= \sum_{\mathbf{k}\in \mathbb{Z}^d} \exp(-2\pi i \mathbf{k}\cdot \bm{\beta}) K_Y(\mathbf{x}-\mathbf{y}+\mathbf{k}\mathbf{m}, \mathbf{0};n),
  \end{equation}
  where $K_Y(\mathbf{x}, \mathbf{y};n)$ denotes the discrete time heat kernel on the graph $Y=\mathcal{C}(\mathbb{Z}^d, S,\pi_S)$.
\end{proposition}

\begin{proof}
  Notice that the sum on the right-hand side of \eqref{eq: heat k on tori as a sum} is finite, and therefore that sum is convergent
  and well-defined. Next, let us check the transformation property under the action of $\mathbf{m}\mathbb{Z}^d$.
  For any $\mathbf{l}\in \mathbb{Z}^d$ we have
\begin{multline*}
   K_{X,\bm{\beta}}(\mathbf{x}+\mathbf{l}\mathbf{m}, \mathbf{y};n)= \sum_{\mathbf{k}\in \mathbb{Z}^d} \exp(-2\pi i (\mathbf{k}+\mathbf{l})\cdot \bm{\beta})   \exp(2\pi i \mathbf{l}\cdot \bm{\beta}) K_Y(\mathbf{x}-\mathbf{y}+(\mathbf{k}+\mathbf{l})\mathbf{m}, \mathbf{0};n)\\= \exp(2\pi i \mathbf{l}\cdot \bm{\beta}) \sum_{\mathbf{k'}\in \mathbb{Z}^d} \exp(-2\pi i \mathbf{k'}\cdot \bm{\beta}) K_Y(\mathbf{x}-\mathbf{y}+\mathbf{k'}\mathbf{m}, \mathbf{0};n)= \exp(2\pi i \mathbf{l}\cdot \bm{\beta})K_{X,\bm{\beta}}(\mathbf{x}, \mathbf{y};n),
  \end{multline*}
  where $\mathbf{k'} = \mathbf{k} + \mathbf{l}$.
The transformation property in the second variable is proved analogously.

It only remains to prove that $  K_{X,\bm{\beta}}(\mathbf{x}, \mathbf{y};n)$ satisfies the diffusion equation \eqref{eq: discrete heat eq inf graph} for the weighted Cayley graph $X=\mathcal{C}(G_{\mathbf{m}},S, \pi_S)$. This follows trivially from the definition of the random walk Laplacian acting on the space of functions initially defined on $X$ and satisfying the transformation property \eqref{eq:twist by character transf prop}.
\end{proof}

Now, we develop the main result of this section, which is a general combinatorial  formula
to compute the heat kernel $K_{X,\bm{\beta}}(\mathbf{x}, \mathbf{0};n)$ of the weighted graph
$X=\mathcal{C}(G_{	\mathbf{m}},\pi _{S},S)$. We assume $S$ is presented as a disjoint union
$S=S_{1}\cup (-S_{1})$, where $-S_{1}$ denotes the set of all the
inverses of the elements of $S_{1}$. We denote by $l>0$ the number of
elements of $S_{1}.$ These elements will be denoted as  $s_{1}=(s_{11},\ldots ,s_{d1}),\ldots ,s_{l}=(s_{1l},\ldots ,s_{dl}).$

In view of Proposition \ref{ heat kernel geometric sum} and Lemma \ref{lem:HK on Zd}, in order to determine $K_{X,\bm{\beta}}(\mathbf{x}, \mathbf{0};n)$ it suffices to count the number of monomials in the variables $t_{1},\ldots t_{d}$ which have the same power mod $\mathbf{m}$  in the polynomial
\[
\left( \sum\limits_{s\in S_{1}} \left(\pi (s)%
\mathbf{t}^{s}+\pi (-s)\mathbf{t}^{-s} \right) \right )
^{n}=\sum\limits_{\substack{(a_{1},\ldots ,a_{l})\in \mathbb{Z}_{\geq 0}^l \\ %
		a_{1}+\cdots +a_{l}=n}}\frac{n!}{a_{1}!\cdots a_{l}!}\prod_{j=1}^{l}\pi(s_j)^{a_j}(\mathbf{t}^{s_1}+\mathbf{t}^{-s_1})^{a_{1}}\cdots (\mathbf{t}^{s_l}+\mathbf{t}^{- s_l})^{a_{l}}.
\]%
Note that we have used the fact that $\pi (s)=\pi (-s)$. This polynomial can be written as
\begin{eqnarray} \notag
	P(t_{1},\ldots ,t_{d}) &=&\sum\limits_{\substack{ \mathbf{a}=(a_{1},\ldots
			,a_{l})\in \mathbb{Z}_{\geq 0}^l  \\ a_{1}+\cdots +a_{l}=n \\ 0\leq j_{r}\leq a_{r}
	}}\frac{n!}{a_{1}!\cdots a_{l}!}\prod_{j=1}^{l}\pi(s_j)^{a_j}\binom{a_{1}}{j_{1}}\cdots
	\binom{a_{l}}{j_{l}}\mathbf{t}^{(2j_1-a_1)s_1}\cdots \mathbf{t}^{(2j_l-a_l)s_l} \notag \\
	&=&\sum\limits_{\substack{ \mathbf{a}=(a_{1},\ldots ,a_{l})\in \mathbb{Z}_{\geq 0}^l  \\ a_{1}+\cdots			+a_{l}=n \\ 0\leq j_{r}\leq a_{r}}}\frac{n!}{%
		a_{1}!\cdots a_{l}!}\prod_{j=1}^{l}\pi(s_j)^{a_j}\binom{a_{1}}{j_{1}}\cdots \binom{a_{l}}{%
		j_{l}}\mathbf{t}_{1}^{2j_{1}-a_{1}}\cdots \mathbf{t}_{l}^{2j_{l}-a_{l}} \label{eq: multivar Poly}
\end{eqnarray}%
where $\mathbf{t}_{1}=t_{1}^{s_{11}}\cdots t_{d}^{s_{d1}},\ldots ,\mathbf{t}%
_{l}=t_{1}^{s_{1l}}\cdots t_{d}^{s_{dl}}.$ Hence, one must determine for any
fixed $l$-tuple $\mathbf{a}=(a_{1},\ldots ,a_{l})\in \mathbb{Z}_{\geq 0}^l $ with $a_1+\ldots+a_l=n$, all possible choices for the vector $(j_{1},\ldots j_{l})$ such that

\begin{equation}
	\begin{pmatrix}
		s_{11} & \cdots  & s_{1l} \\
		\vdots  &  & \vdots  \\
		s_{d1} & \cdots  & s_{dl}%
	\end{pmatrix}%
	\begin{pmatrix}
		2j_{1}-a_{1} \\
		\vdots  \\
		2j_{l}-a_{l}%
	\end{pmatrix}%
	=%
	\begin{pmatrix}
		x_{1}+k_{1}m_{1} \\
		\vdots  \\
		x_{d}+k_{d}m_{d}%
	\end{pmatrix}%
	,  \label{j5}
\end{equation}%
for some $\mathbf{k}=(k_1,\ldots, k_d)\in \mathbb{Z}^d$, under the restrictions $0\leq j_{r}\leq a_{r}.$

Let us denote by $\sigma $ the $\mathbb{Z}$-linear homomorphism $\sigma: \mathbb{Z}^{l}\to\mathbb{Z}^{d} $ determined
by the matrix  $\mathcal{S}=\left(s_{ij}\right)_{d\times l}.$ First, we determine all possible vectors
$z\in \mathbb{Z}^{l}$ such that $p(\sigma (z))=0$, where $p$ denotes the canonical surjection $p: \mathbb{Z}^{d} \to \mathbb{Z}^{d}/\mathbf{m}\mathbb{Z}^{d}=G_{\mathbf{m}}$.

All possible solutions in $G_{\mathbf{m}}$ of this problem can be
obtained as follows. By selecting appropriate bases $B_{\mathbb{Z}^{l}}=\{v_{i}:i=1,\ldots ,l\}$ and $B_{\mathbb{Z}^{d}}=\{w_{j}:j=1,\ldots
,d\}$ for $\mathbb{Z}^{l}$ and $\mathbb{Z}^{d}$, respectively, the matrix $\mathcal{S}$ has diagonal form when written in these bases
\[
(\mathcal{S} )=%
\begin{pmatrix}
	q_{1} & \cdots  & 0 & \cdots  \\
	\vdots  & \ddots  & \vdots  & \vdots  \\
	0 & \cdots  & q_{t} & \cdots
\end{pmatrix}_{t\times l}%
\]%
Here $t=\max\{d,l\}$ and $\{q_{i}\}$ denotes the set of invariant factors of $\sigma$. In case $d<l$, we define $q_{d+1}=\ldots =q_{l}=0$. Finding these basis amounts to finding invertible matrices $U$ and $V$ such that $U\mathcal{S} V= (\mathcal{S} )$ is the Smith normal form of $\mathcal{S} $; see \cite{Sm61}.

Hence, for $i=1,\ldots,l$ one has $\sigma (v_{i})=q_{i}w_{i}.$ Using the matrix $U=(u_{ij})_{d\times d}$ we can
write $w_{j}=u_{1j}e_{1}+\cdots +u_{dj}e_{d}$, where $e_{1},\ldots ,e_{d}$
denotes the canonical basis for $\mathbb{Z}^{d}.$ Let $(z_1,\ldots,z_l)$ be the coordinates of $z\in\mathbb{Z}^l$ in the basis $\{v_1,\ldots,v_l\}$. Solving for $p(\sigma(z))=p(\sigma (\sum\limits_{i=1}^{l}z_{i}v_{i}))=0$
is equivalent to solving in $G_{\mathbf{m}}$ the system of $d$ equations
\begin{equation}
	\begin{bmatrix}
		z_{1}c_1\equiv 0\text{ mod }m_{1} \\
		\vdots  \\
		z_{d}c_d\equiv0\text{ mod }m_{d}%
	\end{bmatrix}%
	,  \label{j8}
\end{equation}%
%(Note that in case $d < i\leq l$ we have $q_{i}=0,$ by definition, hence all sums in \eqref{j8} reduce to first $d$ terms.)
where $c_{i}=\sum\limits_{j=1}^{d}(u_{ij}q_{i})$, $i=1,\ldots,d$. Each equation $
c_{i}z_{i}\equiv0$ mod $m_{i}$, $i=1,\ldots,d$ in the system \eqref{j8} can be readily solved in $
\mathbb{Z}/m_i\mathbb{Z}$ by elementary number theory. Thus, we assume we have already computed all possible solutions in $G_{	\mathbf{m}}$ of (\ref{j8}) and they are given as the $N$ elements of the set
\[
L_{S,\mathbf{m}}=\left\{ \mathbf{z}_{\mu }= (z_{1\mu}, \ldots, z_{d\mu}):\mu =1,\ldots ,N\right\} .
\]
Now, in order to solve (\ref{j5}) we can take each solution $(z_{1\mu}, \ldots, z_{d\mu})\in L_{S,\mathbf{m}}$ and find all possible $\mathbf{k=(}k_{1},\ldots
,k_{d})$ such that
\[
\begin{array}{c}
	2j_{1}-a_{1}=z_{1\mu }+x_1+k_{1}m_{1} \\
	\vdots  \\
	2j_{d}-a_{d}=z_{d\mu }+x_d+k_{d}m_{d}%
\end{array}%
,
\]

where the original restrictions $0\leq j_{r}\leq a_{r}$ are equivalent to
choosing $\mathbf{k}$ such that
\begin{eqnarray}
	0 \leq a_{r}+z_{r\mu }+x_r+k_{r}m_{r}\leq  2a_{r},\text{ for }r=1,\ldots ,d,
	\label{r1} \\
	0 \leq j_{r}\leq a_{r}\text{, for }r=d+1,\ldots ,l\text{, in case }d<l.
	\label{r2}
\end{eqnarray}%
Thus, $\mathbf{k}$ must satisfy for each $\mu =1,\ldots,N$ that
\[
\left\vert z_{r \mu} +x_r+k_{r}m_{r}\right\vert \leq a_{r}, \text{ for }r=1,\ldots
,d.
\]%

With all this discussion, we have proved the following theorem.

\begin{theorem}\label{thm:hk_formula_quotient}  With notation as above, let
$$
A_{ \bm{\beta},S,\mathbf{m}}(a_1,\ldots, a_l;n; \mathbf{k})= e^{-2\pi i
\mathbf{k\cdot \bm{\beta} }}\frac{n!}{a_{1}!\cdots a_{l}!}\prod_{j=1}^{l}\pi(s_j)^{a_j} \prod\limits_{r=1}^{d}\binom{a_{r}}{(a_{r}+x_{r}+z_{r\mu }+k_{r}m_{r})/2}.
$$
Then the heat kernel $K_{G_{\mathbf{m}},S,\bm{\beta }}(\mathbf{x},0;n)$ is given by the following formulas.
\begin{itemize}
\item[1.] If $d\geq l$:
$$\sum\limits_{_{\substack{ (a_{1},\ldots ,a_{l})\in \mathbb{Z}_{\geq 0}^l \\ a_{1}+\cdots +a_{l}=n}}}\sum\limits_{\mathbf{z}_{\mu }\in L_{S,\mathbf{m}%
}}\sum\limits_{\substack{ _{\substack{ \mathbf{k}\in \mathbb{Z}^{d} \\ %
\left\vert z_{r\mu }+x_{r}+k_{r}m_{r}\right\vert \leq a_{r},\text{ }%
r=1,\ldots ,d}} \\ a_{r}+x_{r}+z_{r\mu }+k_{r}m_{r}\text{ even }}}A_{ \bm{\beta},S,\mathbf{m}}(a_1,\ldots, a_l;n; \mathbf{k}),$$
\item[2.] If $d<l$
$$\sum\limits_{_{_{\substack{ (a_{1},\ldots ,a_{l})\in \mathbb{Z}_{\geq 0}^l \\ a_{1}+\cdots +a_{l}=n \\ 0\leq \text{ }j_{r}\text{ }\leq a_{r},\
d+1\leq r\leq l}}}}\sum\limits_{\mathbf{z}_{\mu }\in L_{S,\mathbf{m}%
}}\sum\limits_{\substack{ _{\substack{ \mathbf{k}\in \mathbb{Z}^{d} \\ %
\left\vert z_{r\mu }+x_r+k_{r}m_{r}\right\vert \leq a_{r},\text{ }r=1,\ldots ,d
\\ 0\leq j_{r}\leq a_{r}\text{, }r=d+1,\ldots ,l\text{ }}} \\ %
a_{r}+x_{r}+z_{r\mu }+k_{r}m_{r}\text{ even, }r=1,\ldots ,d,\text{ }}}
\prod\limits_{r=d+1}^{l}\binom{a_{r}}{j_{r}}A_{ \bm{\beta},S,\mathbf{m}}(a_1,\ldots, a_l;n; \mathbf{k}).
$$
\end{itemize}
\end{theorem}
\begin{proof}
	It readily follows from (\ref{eq: heat k on tori as a sum}), (\ref{r1}) and (\ref{r2}).
\end{proof}

\begin{example}\rm \label{ex. d=1 whole group HK}
Let $d=1$, $S=\{\pm 1\}$, $\pi_S(s)=1/|S| =1/2$ for $s\in S$. Then an application
of Lemma \ref{lem:HK on Zd} and the binomial theorem yields that the heat kernel $K_{Y}(x,0;n)$ for $n\geq 0$ on the Cayley graph  $Y=\mathcal{C}(\mathbb{Z}, S,\pi_S)$ is given by
\begin{equation}\label{eq: heat k on Z - formula}
K_{Y}(x,0;n)= \left\{
                       \begin{array}{ll}
                         2^{-n}\binom{n}{(x+n)/2}, & \text{  if  } |x|\leq n \text{  and  } x+n \text{  is even}; \\
                         0, & \text{  if  } |x|>n \text{  or  } x+n \text{  odd}.
                       \end{array}
                     \right.
\end{equation}
This is a well-known result from probability theory; see \cite{Wo00} as well as \cite{SS14}, Example 3.3 for a related, more general result.
\end{example}

\begin{example}\rm \label{ex: d=1, b any beta any}
Let $m$ be a positive integer, and let $b\in\{1,\ldots, m-1\}$. Consider the weighted Cayley graph $X=\mathcal{C}(\mathbb{Z}/m\mathbb{Z}, S, \pi_S)$,
where $S=\{-b,b\}$ is the symmetric subset of $\mathbb{Z}/m\mathbb{Z}$ and $\pi_S(s)=1/2$, for $s\in S$. Let $Y$ be the weighted Cayley graph $\mathcal{C}(\mathbb{Z}, S, \pi_S)$. Then, from Proposition \ref{ heat kernel geometric sum} with $d=l=1$ and arbitrary
$\beta\in\mathbb{R}$, the discrete time heat kernel on $X$ twisted by the character $\chi_{\beta}$ on $\mathbb{Z}$ is given by
\begin{equation}\label{eq: K of example}
K_{X,\beta}(x,y;n)=\sum_{\ell=-\infty}^{\infty} e^{-2\pi i \ell \beta} K_Y(x-y+\ell m,0;n),
\end{equation}
where
\begin{equation}\label{eq: KY of example}
 K_Y(z+\ell m,0;n)=\left\{
              \begin{array}{ll}
                2^{-n}\binom{n}{(z+m\ell + bn)/2b}, & \text{  if  } |z+m\ell|\leq bn \text{  and  } 2b \mid(z+m\ell +bn) \\ \\
                         0 &

                          \text{  otherwise.}
              \end{array}
            \right.
\end{equation}
\end{example}

\begin{remark}\rm
 An application of Theorem \ref{thm:hk_formula_quotient} with $d=1$, $\bm{\beta}=0$ and $S=\{\pm1, \ldots, \pm m\}$ yields
 the combinatorial expression for the discrete time and discrete space heat kernel $u(x,t)$ for $x\in \mathbb{Z}$ and $t\in \mathbb{Z}_\geq 0$,
 which is a solution to the diffusion equation
 $$
 u(x,t+1)-u(x,t)=\sum_{i=-m}^m a_i u(x+i,t)\quad \text{with}\quad\sum_{i=-m}^m a_i=0,
 \text{  } a_0\leq 0, \text{  } a_i\geq 0 \text{\,\, for\,\,}  i\neq 0
 $$
 as studied in \cite{FSS14}.  Moreover, Theorem \ref{thm:hk_formula_quotient}
 provides a means by which one can explicitly evaluate solutions to the diffusion equation on $\mathbb{Z}^d$ for any $d>1$.
 %Several applications of this formula will be given in the subsequent paper \cite{CHJSV22a}.
 \end{remark}

\section{The main result}

In this section we present our main result, which is equation \eqref{eq: main formula}.
The formula relates a sum of products of binomial coefficients to a sum of $n$-th powers of trigonometric functions,
after each side of the equation is twisted by a character.
The proof is based on the existence of two different expressions for the discrete time twisted heat kernel: First, from a purely combinatorial
expression, given in Proposition  \ref{ heat kernel geometric sum} and Theorem \ref{thm:hk_formula_quotient}; and
second, using spectral theory, stemming from the fact that the
eigenvalues of the adjacency matrix are expressed in terms of the characters of the underlying group of the corresponding Cayley
graph as proved in Corollary 3.2 of \cite{Ba79}.  Because the heat kernel is unique, any two distinct formulas for the
heat kernel must be equal.

%The spectrum of the matrix $A_X$ can be expressed in terms of a linear combination of products of characters on the groups $\mathbb{Z}/m_j\mathbb{Z}$, $j=1,\ldots, d$. Thus, the spectral expression of the heat kernel at discrete time $n$ can be written as a sum of products of shifted $n$th powers of cosine functions (due to symmetry of $S$).

%On the other hand, the arithmetic expression can be written as a linear combination of certain products of binomial coefficients (stemming from the twisted by $\chi_{\bm{\beta}}$ heat kernel on the covering space $\mathbb{Z}^d$). The uniqueness of the heat kernel yields that the two expressions must be equal.

%\textcolor{red}{We should explain this better since this is the heart of our results.}

\begin{theorem} \label{thm: main} Let $d\geq 1$ be an integer, $\mathbf{m}=(m_1,\ldots,m_d)\in\mathbb{Z}_{>0}^d$
and $\bm{\beta}=(\beta_1,\ldots , \beta_d)\in \mathbb{R}^d$. Let $G_{\mathbf{m}}$ denote the discrete torus
$\mathbb{Z}^d/\mathbf{m}\mathbb{Z}^d$ with $M=\prod_{j=1}^d m_j$ vertices, and let $S$ be a symmetric subset of
$G_{\mathbf{m}}$ with an associated positive probability distribution $\pi_S$. Then, for any positive integer
$n$ and $\mathbf{x},\mathbf{y}\in \mathbb{Z}^d$, we have that
 \begin{equation} \label{eq: main formula}
  \sum_{\mathbf{k}\in \mathbb{Z}^d} \exp(-2\pi i \mathbf{k}\cdot \bm{\beta}) K_Y(\mathbf{x}-\mathbf{y}+\mathbf{k}\mathbf{m}, \mathbf{0};n) = \frac{1}{M} \sum_{\mathbf{n}\in G_{\mathbf{m}}}(\lambda_{\mathbf{n}}^{\mathbf{m}, \bm{\beta}})^n \chi_{\mathbf{n}}^{\mathbf{m}, \bm{\beta}}(\mathbf{x})\overline{\chi_{\mathbf{n}}^{\mathbf{m}, \bm{\beta}}(\mathbf{y})},
  \end{equation}
  where $K_Y$ is defined in Lemma \ref{lem:HK on Zd}, $\{\chi_{\mathbf{n}}^{\mathbf{m}, \bm{\beta}}\}_{\mathbf{n}=(n_1,\ldots, n_d) \in G_{\mathbf{m}}}$ is the set of additive characters on $\mathbb{Z}^d$ defined by
   \begin{equation}\label{eq: defn of chi,n}
\chi_{\mathbf{n}}^{\mathbf{m}, \bm{\beta}}(\mathbf{x})=\prod_{j=1}^d \exp\left(2\pi i \frac{n_j+\beta_j}{m_j} x_j\right), \quad \text{for all} \quad \mathbf{x}=(x_1,\ldots,x_d) \in \mathbb{Z}^d,
\end{equation}
      and $\lambda_{\mathbf{n}}^{\mathbf{m}, \bm{\beta}}$, for $\mathbf{n}\in G_{\mathbf{m}}$ is given by
   \begin{equation}\label{eq: defn of lambda,n}
\lambda_{\mathbf{n}}^{\mathbf{m}, \bm{\beta}}=\sum_{s\in S} \pi_S(s) \chi_{\mathbf{n}}^{\mathbf{m}, \bm{\beta}}(s).
\end{equation}
\end{theorem}

\begin{proof}

The proof follows from the spectral expansion of the discrete time twisted heat kernel. Namely, the eigenfunctions of the
adjacency matrix of the Cayley graph $X=\mathcal{C}(G_{\mathbf{m}},S,\pi_S)$ correspond to the characters of the abelian
group $G_{\mathbf{m}}$. Therefore, it is straightforward that $\{\frac{1}{\sqrt{M}}\chi_{\mathbf{n}}^{\mathbf{m}, \bm{\beta}}\}_{\mathbf{n}\in
G_{\mathbf{m}}}$ is a set of $L^{2}$-normalized eigenfunctions of the adjacency matrix  $A_X$  which
satisfy the transformation property \eqref{eq:twist by character transf prop}.
The corresponding  eigenvalues are $\{\lambda_{\mathbf{n}}^{\mathbf{m}, \bm{\beta}}\}_{\mathbf{n} \in G_{\mathbf{m}}}$ ; there are $M$ of them, counted with multiplicities.

The matrix which corresponds to the random walk on the weighted graph $X$ acting on the set of functions
satisfying the transformation property \eqref{eq:twist by character transf prop} is diagonalized
when using the eigenfunctions $\frac{1}{\sqrt{M}} \chi_{\mathbf{n}}^{\mathbf{m}, \bm{\beta}}$ as a basis. Hence, we
may write the discrete time heat kernel twisted by $\chi_{\bm{\beta}}$ at time $n\geq 1$ as
\begin{equation}\label{eq: heat k spectral exp}
  K_{X,\bm{\beta}}(\mathbf{x}, \mathbf{y};n)=\frac{1}{M} \sum_{\mathbf{n}\in G_{\mathbf{m}}}(\lambda_{\mathbf{n}}^{\mathbf{m}, \bm{\beta}})^n \chi_{\mathbf{n}}^{\mathbf{m}, \bm{\beta}}(\mathbf{x})\overline{\chi_{\mathbf{n}}^{\mathbf{m}, \bm{\beta}}(\mathbf{y})}.
\end{equation}
Thus, equation \eqref{eq: main formula} follows from the uniqueness of the heat kernel, once we combine
equation \eqref{eq: heat k spectral exp} and Proposition \ref{ heat kernel geometric sum}.
\end{proof}

A particularly simple combinatorial identity can be deduced by taking
$\mathbf{x}= \mathbf{y}$ in \eqref{eq: main formula}.  Indeed, each
character takes values in the set of complex numbers of absolute value equal to one,
so then
$$
 \chi_{\mathbf{n}}^{\mathbf{m}, \bm{\beta}}(\mathbf{x})\overline{\chi_{\mathbf{n}}^{\mathbf{m}, \bm{\beta}}(\mathbf{x})}=1,
$$
for all $\mathbf{x}$, $\mathbf{n}$, $\mathbf{m}$ and $\bm{\beta}$.  Therefore, we get the following corollary.

\begin{corollary}\label{cor: main trace f-la}
With notation as above and any positive integer $n$, we have that
    \begin{equation} \label{eq: main formula - cor}
\sum_{\mathbf{n}\in G_{\mathbf{m}}}(\lambda_{\mathbf{n}}^{\mathbf{m}, \bm{\beta}})^n=  M \sum_{\mathbf{k}\in \mathbb{Z}^d} \exp(-2\pi i \mathbf{k}\cdot \bm{\beta}) K_Y(\mathbf{k}\mathbf{m}, \mathbf{0};n),
  \end{equation}
  where $K_Y$ is defined in Lemma \ref{lem:HK on Zd} and  $\lambda_{\mathbf{n}}^{\mathbf{m}, \bm{\beta}}$ is given by \eqref{eq: defn of lambda,n}.
\end{corollary}

Formula \eqref{eq: main formula - cor} shows how to compute explicitly the trigonometric sum which is the trace of the discrete time heat kernel of the random walk Laplacian, twisted by $\chi_{\bm{\beta}},$ in terms of a twisted sum of binomial coefficients. Further simplifications of formula \eqref{eq: main formula - cor} can be obtained by taking $\bm{\beta} \in \mathbf{m}\mathbb{Z}^d$ or $\bm{\beta} \in \frac{1}{2} \mathbf{m}\mathbb{Z}^d$. Namely, when $\bm{\beta} \in \mathbf{m}\mathbb{Z}^d$, we have $\exp(-2\pi i \mathbf{k}\cdot \bm{\beta})=1$, while for  $\bm{\beta} \in \frac{1}{2} \mathbf{m}\mathbb{Z}^d$ we have $\exp(-2\pi i \mathbf{k}\cdot \bm{\beta})=(-1)^{\mathbf{k}\cdot \mathbf{m}}$. The expressions for $\lambda_{\mathbf{n}}^{\mathbf{m}, \bm{\beta}}$ can be significantly simplified in those cases as well.

\begin{remark}\label{rem:Galois_argument}\rm
Assume $\pi _{S}$ is the uniform distribution and that there is no
character twist, meaning $\bm{\beta }=0$.  Then the left hand side of equation \eqref{eq: main formula - cor}
must be an integer. To see this, let $E$ denote the normal closure over $\mathbb{Q}$ containing
all the primitive $m_{i}$-th roots of unity $\omega _{i}$ for $i=1,\ldots ,d$.  Any
element $\phi $ of the Galois group Gal$(E/\mathbb{Q})$ maps any $%
\omega _{i}$ into another primitive $m_{i}$-th root of unity, meaning $\phi
(\omega _{i})=\omega _{i}^{r_{i}}$ where $0<r_{i}<m_{i}$ is relatively prime
to $m_{i}$. Therefore,
\begin{eqnarray*}
\phi \left(|S| \lambda _{\mathbf{n}}^{\mathbf{m},\mathbf{0}}\right)
&=&\sum\limits_{s=(s_{1},\ldots ,s_{d})\in S}\phi \left( \chi _{\mathbf{n}}^{\mathbf{m,0}%
}\right) =\sum\limits_{(s_{1},\ldots ,s_{d})\mathbf{\in }S}\phi (\omega
_{1}^{n_{1}s_{1}}\cdots \omega _{d}^{n_{d}s_{d}}) \\
&=&\sum\limits_{(s_{1},\ldots ,s_{d})\mathbf{\in }S}\omega
_{1}^{r_{1}n_{1}s_{1}}\cdots \omega
_{d}^{r_{d}n_{d}s_{d}}=\sum\limits_{(s_{1},\ldots ,s_{d})\mathbf{\in }%
S}\omega _{1}^{n_{1}^{\prime }s_{1}}\cdots \omega _{d}^{n_{d}^{\prime }s_{d}}=|S|\lambda _{\mathbf{n}^{\prime }}^{\mathbf{m},\mathbf{0}},
\end{eqnarray*}%
where $|S|$ is the cardinality of $S$, $\mathbf{n}^{\prime }=\mathbf{rn},$ and $\mathbf{r}=(r_{1},\ldots
,r_{d}).$ Since each $r_{i}$ is relatively prime to $m_{i}$,
multiplication by $\mathbf{r}$ permutes the elements of $\mathbb{Z}_{\mathbf{%
m}}$.  Consequently, $\phi $ permutes the set of eigenvalues $\lambda _{%
\mathbf{n}}^{\mathbf{m},\mathbf{0}}$. Hence, if $a=|S|^n\sum\limits_{\mathbf{%
n}\in G_{\mathbf{m}}}(\lambda _{\mathbf{n}}^{\mathbf{m},\mathbf{0}})^{n}$ then
\[
\phi (a)=|S|^n\sum\limits_{\mathbf{n}^{\prime }\in G_{\mathbf{m}}}(\lambda _{%
\mathbf{n}^{\prime }}^{\mathbf{m},\mathbf{0}})^{n}=a,
\]%
and, therefore, $a$ must be a rational number. Clearly $a$ is an algebraic
integer, so then $a$ must be an integer. This shows that one interesting consequence of
our main theorem is that it provides an explicit formula to compute the integer $a$ in terms of binomial
coefficients.

Implicit in this argument is the fact that the left-hand side of \eqref{eq: main formula - cor} is a symmetric polynomial in the eigenvalues $\lambda_{1},\ldots,\lambda_{M}$ defined in \eqref{eq: defn of lambda,n}, where $M$ denotes the order of $G_{\bold{m}}.$ In this case the symmetric polynomial is $p_n(x_{1},\ldots,x_{M})=x_{1}^n+\cdots+x_{M}^n.$  It is a well known fact that any symmetric polynomial in the variables $x_{1},\ldots,x_{M}$ can be written as a polynomial with rational coefficients in the $p_{k}$'s, where $p_k(x_{1},\ldots,x_{M})=x_{1}^k+\cdots+x_{M}^k$, with $k=1,\ldots,n$ \cite{STF99}. Therefore, for any symmetric polynomial with rational coefficents $q(x_{1},\ldots,x_{M})$, the value of $q(\lambda_{1},\ldots,\lambda_{M})$ \emph{is a rational number}.
% An explicit algorithm to compute this rational number will be developed in the subsequent article \cite{CHJSV22a}.
\end{remark}

\section{Trigonometric sums stemming from the twisted heat kernel on a discrete torus}\label{sec: trig sums}

In this section, we present corollaries of our main result that can be deduced from \eqref{eq: main formula} by an appropriate choice of $d,$ $\mathbf{m}$ and $\bm{\beta}$. We derive many interesting identities, and show that special instances of those identities are main results from e.g. \cite{dFGK17}, \cite{dFK13}, \cite{Me12} and \cite{BH10}.

\subsection{Trigonometric sums twisted by an additive character}

The first corollary in this section is a consequence of the computation of the trace of the heat kernel, as in Corollary \ref{cor: main trace f-la}.

\begin{corollary} \label{ex: d=1}
For arbitrary positive integers $m,n$ and arbitrary $\beta\in \mathbb{R}$,
\begin{equation} \label{eq: main identity d=1}
\sum\limits_{j=0}^{m-1}\cos ^{n} \left(\frac{2\pi (j+\beta)}{m}\right)=2^{-n}m\sum\limits_{k=-%
\left\lfloor n/m\right\rfloor ,\text{ \ }km+n\text{ even}}^{\left\lfloor
n/m\right\rfloor }e^{-2\pi i k\beta} \binom{n}{(km+n)/2}.
\end{equation}
Here $\lfloor x\rfloor$ stands for the largest integer less than or equal to $x$.
\end{corollary}
\begin{proof}
 We apply Corollary \ref{cor: main trace f-la} with $d=1$, $\mathbf{m}=m=M$, $\bm{\beta}=\beta$, $S=\{-1,1\}$ and $\pi_S(s)=1/2$, for $s\in S$. Then, for a fixed, arbitrary integer $n>0$, an application of formula \eqref{eq: heat k on Z - formula} with $x=km$ shows that, for any
% gives
%$$
%K_Y(km,0;n)=\left\{
%              \begin{array}{ll}
%                2^{-n}\binom{n}{(km+n)/2}, & \text{  if  } |km|\leq n \text{  and  } km+n %\text{  is even}; \\
%                         0, & \text{  if  } |km|>n \text{  or  } km+n \text{  odd}.
%              \end{array}
%            \right.
%$$
$\beta\in \mathbb{R}$, the right hand side of \eqref{eq: main formula - cor} equals
$$
2^{-n}m\sum\limits_{k=-
\left\lfloor n/m\right\rfloor ,\text{ \ }km+n\text{ even}}^{\left\lfloor
n/m\right\rfloor }e^{-2\pi i k\beta} \binom{n}{(km+n)/2}.
$$
For $j\in 0,\ldots,m$ the eigenvalues $\lambda_j^{m,\beta}$ of the adjacency matrix are
$$
\lambda_j^{m,\beta} = \frac{1}{2}\left(e^{2\pi i (j+\beta)/m} +e^{-2\pi i (j+\beta)/m} \right)=\cos(2\pi  (j+\beta)/m).
$$
Thus, formula \eqref{eq: main formula - cor} yields identity \eqref{eq: main identity d=1}.
\end{proof}

Next, we examine a few special choices of $\beta$, in order to show that \eqref{eq: main identity d=1} generalizes the
main results in \cite{dFGK17} and \cite{Me12}.

\noindent
1. For $\beta=0$, specializing \eqref{eq: main identity d=1} to the case where $m=2m_{1}$ and $n=2n_{1}$ are even natural numbers we get
\begin{equation*}
\sum\limits_{j=0}^{2m_{1}-1}\cos ^{2n_{1}} \left(\frac{\pi j}{m_{1}}%
\right)=2^{-2n_{1}}2m_{1}\sum\limits_{k=-\left\lfloor n/m\right\rfloor ,\text{ }%
}^{\left\lfloor n/m\right\rfloor }\binom{2n_{1}}{km_{1}+n_{1}},
\end{equation*}
which can be simplified to give equation \eqref{eq: example1} which is the first main result of \cite{dFGK17}.

\noindent
2. When $\beta=-m/4$ , we have $\cos \left(\frac{2\pi (j+\beta)}{m}\right)= \sin (2\pi j/m)$, and \eqref{eq: main identity d=1} becomes
\begin{equation*}
\sum\limits_{j=0}^{m-1}\sin ^{n}\left(\frac{2\pi j}{m}\right)=2^{-n}m\sum\limits_{k=-%
\left\lfloor n/m\right\rfloor ,\text{ \ }km+n\text{ even}}^{\left\lfloor
n/m\right\rfloor }i^{km}\binom{n}{(km+n)/2}.
\end{equation*}
By specializing to the case where $m=2m_{1}$ and $n=2n_{1}$ are even natural
numbers, analogously as above we obtain
\begin{equation*}
\sum\limits_{j=0}^{m_{1}-1}\sin ^{2n_{1}}\left(\frac{\pi j}{m_{1}}%
\right)=2^{-2n_{1}}m_{1}\sum\limits_{k=-\left\lfloor n_{1}/m_{1}\right\rfloor ,%
\text{ }}^{\left\lfloor n_{1}/m_{1}\right\rfloor } (-1)^{km_1}\binom{2n_{1}}{km_{1}+n_{1}%
}.
\end{equation*}
This proves the second main result of \cite{dFGK17}.

\noindent
3. By choosing $\beta=0$, $m=2m_{1}$ and $n=2p$ even natural numbers, and applying trigonometric identities to reduce the range of summation, one easily deduces Theorem 1 from \cite{Me12}. Theorem 3 from \cite{Me12} follows by the same procedure, after taking $m=2m_{1}$ and $n=2p$ even natural numbers and $\beta=1/2$ in \eqref{eq: main identity d=1}.
%\item Describe what to chose to get the main result of \cite{Ho18}...

In the following corollary we apply the main theorem for $d=1$ and show how to deduce a combinatorial expression for powers of sines and cosines twisted by an additive character.

\begin{corollary} \label{ex: d=1, with character}
For arbitrary positive integers $m,n$, any integer $b\in\{1,\ldots,m-1\}$, arbitrary $\alpha,\varphi\in \mathbb{R}$ and any character $\chi$ of the additive group $\mathbb{Z}/m\mathbb{Z}$, we have that
\begin{equation}\label{eq: cos twisted sum}
\sum_{j=0}^{m-1}\chi(j) \cos^n \left(\frac{2\pi j b}{m} + \alpha\right)=\frac{m}{2^n} \sum\limits_{_{\substack{ d \in\{0,\ldots, n\} \\ m\mid (2d-n)b-r }}}\binom{n}{d} \exp\left(i\alpha (n-2d)\right),
\end{equation}
and
\begin{equation}\label{eq: sin twisted sum}
\sum_{j=0}^{m-1}\chi(j) \sin^n \left(\frac{2\pi j b}{m} + \varphi \right)=\frac{m}{2^n} \sum\limits_{_{\substack{ d \in\{0,\ldots, n\} \\ m\mid (2d-n)b-r }}}\binom{n}{d} \exp\left(i\varphi(n-2d)\right)\exp\left(-i\frac{\pi}{2}(n-2d)\right).
\end{equation}
Here $r \in\{1,\ldots,m-1\}$ is an integer such that  $\chi(x)=\exp(2\pi i rx/m)$ for all $x\in \mathbb{Z}/m\mathbb{Z}$ and the empty sum equals zero.
\end{corollary}

\begin{proof}
We use the setting of Example \ref{ex: d=1, b any beta any} in which the discrete time heat kernel on $X$ twisted by the character $\chi_{\beta}$ on $\mathbb{Z}$ is given for positive integers $n$ and $x,y \in \mathbb{Z}/m\mathbb{Z}$ by \eqref{eq: K of example} and \eqref{eq: KY of example}. Therefore, the left hand side of formula \eqref{eq: main formula} becomes:
\begin{equation} \label{eq: LHS with character}
K_{X,\beta}(x,y;n)=  \sum\limits_{_{\substack{ d \in\{0,\ldots, n\} \\ m\mid (2d-n)b-(x-y) }}}2^{-n} \binom{n}{d} \exp(2\pi i (x-y + (n-2d)b)\beta/m).
\end{equation}
For $j\in\{0,\ldots,m\}$ the product of twisted characters $\chi_j^{m,\beta}$ on the right-hand side of formula \eqref{eq: main formula} is
$$
\chi_j^{m,\beta}(x)\overline{\chi_j^{m,\beta}(y)}=\exp(2\pi i (j+\beta)(x-y)/m)
$$
and the associated eigenvalues are
$$
\lambda_j^{m,\beta}= \frac{1}{2}\left(\exp(2\pi i (j+\beta)b/m) + \exp(-2\pi i (j+\beta)b/m) \right) = \cos \left(\frac{2\pi j b}{m} + \frac{2\pi \beta b}{m}\right).
$$
Equation \eqref{eq: main formula} now yields the following identity, valid for any real number $\beta$,  positive integers $n$ and $x,y \in \mathbb{Z}/m\mathbb{Z}$:
$$
\sum_{j=0}^{m-1}\exp\left(2\pi i \frac{j(x-y)}{m}\right) \cos^n \left(\frac{2\pi j b}{m} + \frac{2\pi \beta b}{m}\right)=\frac{m}{2^n} \sum\limits_{_{\substack{ d \in\{0,\ldots, n\} \\ m\mid (2d-n)b-(x-y) }}}\binom{n}{d} \exp\left(2\pi i \frac{(n-2d)b\beta}{m}\right).
$$
We take $x,y$ above so that $r=x-y$, and deduce, by choosing $\beta=\frac{m \alpha }{2\pi b}$, the identity \eqref{eq: cos twisted sum}.

By choosing $\alpha=\varphi-\frac{\pi}{2}$, equation \eqref{eq: cos twisted sum} becomes  \eqref{eq: sin twisted sum}. This completes the proof of the corollary.
\end{proof}

The following example illustrates the results which are stated in Corollary \ref{ex: d=1, with character}.

\begin{example}\rm
For some positive integer $k$, let
$m=3k$ and $r=k$, and set $\alpha=0$ in Corollary \ref{ex: d=1, with character}. Then, one gets the following interesting formula for the sum of powers of cosines twisted by the third root of unity $\omega=\exp(2\pi i/3)$:
$$
\sum_{j=0}^{3k-1} \omega^j\cos^n \left(\frac{2\pi j }{3k}\right)=\frac{3k}{2^n} \sum\limits_{_{\substack{ d \in\{0,\ldots, n\} \\ 3k\mid (2d-n-k) }}}\binom{n}{d}.
$$
We find it fascinating that the sum on the left hand side is a rational number. For example, when $m=102$, so $k=34$, and
taking $n=100$ one has that
$$
\sum_{j=0}^{101} \omega^j\cos^{100} \left(\frac{\pi j }{51}\right)=\frac{102}{2^{100}}\left(\binom{100}{16} + \binom{100}{67}\right)=\frac{7514656923394238847040235025}{316912650057057350347175801344}.
$$
\end{example}

Formulas \eqref{eq: cos twisted sum} and \eqref{eq: sin twisted sum} can be combined with simple trigonometric identities needed to reduce the range of summation to yield the new, succinct expression for the evaluation of the alternating sum of powers of cosines with fractional multiples of $\pi/2$ evaluated in \cite{dFK13}. Namely, we have the following corollary.

\begin{corollary}\label{cor:NewNumerics}
For two positive integers $m,n$, let $$S(n,m):= \sum\limits_{k=1}^m (-1)^k \cos^{2n}\left(\frac{k\pi}{2m+2}\right).$$
Then
\begin{equation} \label{eq: halp pi sum of cos}
S(n,m)=-\frac{1}{2}+\frac{m+1}{4^n}\sum\limits_{_{\substack{ d \in\{0,\ldots, 2n\} \\ 2(m+1)\mid d-n-(m+1) }}}\binom{2n}{d},
\end{equation}
with the convention that the empty sum equals zero.
\end{corollary}
\begin{proof}
  Let $M:=m+1$. A simple change of variables in the summation gives that
  $$
  \sum\limits_{k=0}^{4M-1} (-1)^k \cos^{2n}\left(\frac{k\pi}{2M}\right)=2 + 2S(n,m)+2 \sum\limits_{k=0}^{M-1} (-1)^{k+M} \sin^{2n}\left(\frac{k\pi}{2M}\right)=2 + 4S(n,m)
  $$
  where the second equality follows by expressing $S(n,m)$ in terms of sums of sines (as in \cite{dFK13}, page 357).
On the other hand, equation \eqref{eq: cos twisted sum} with $m=4M$, $b=1$, $\alpha=0$, $n=2n$ and $r=2M$ yields
$$
\sum\limits_{k=0}^{4M-1} (-1)^k \cos^{2n}\left(\frac{k\pi}{2M}\right)= \frac{4M}{4^n}\sum\limits_{_{\substack{ d \in\{0,\ldots, 2n\} \\ 4M\mid 2d-2n-2M }}}\binom{2n}{d}.
$$
Combining the last two displays yields identity \eqref{eq: halp pi sum of cos}.
\end{proof}

%Proceeding in a similar manner as above one can deduce Theorem 2 of \cite{Me12}.

%\vskip .06in
Let us note that the procedure to compute $S(n,m)$ using the identity \eqref{eq: halp pi sum of cos} is straightforward.
For example, when $n=100$ and $m=13$, it suffices to determine the set of $d\in\{0,1,\ldots,200\}$ for which $28\mid (d-114)$; this set equals $\{2,30,58,86,114,142,170,198\}$. Hence,
$$
S(100,13)=-\frac{1}{2}+ \frac{14}{4^{100}}\cdot 2\left(\binom{200}{2}+\binom{200}{30} + \binom{200}{58}+\binom{200}{86}\right).
$$
Inserting this expression into the \emph{Wolfram Alpha} online calculator, and after a few seconds, one obtains
$$
S(100,13)= - \frac{27820144416504768614943952313424972252178190684153272739503}{100433627766186892221372630771322662657637687111424552206336}
$$
which exactly matches the evaluation of the same sum on p. 373 of \cite{dFK13}.
As a further example (not computed in \cite{dFGK17}) we get that
$$
S(110,18)=-\frac{1}{2}+ \frac{19}{4^{110}}\cdot 2\left(\binom{220}{15}+\binom{220}{53} + \binom{220}{91}\right).
$$
The \emph{Wolfram Alpha} calculator produces the rational number
$$
S(110,18)=-\frac{89182224882179103045185472064334917993187398846393647267026811637}{210624583337114373395836055367340864637790190801098222508621955072}.
$$

\subsection{Trigonometric sums twisted by a multiplicative character}

Equation \eqref{eq: sin twisted sum} yields the main result of \cite{BH10}, their Theorem 1.2,
which is also related to \cite{BZ04} and \cite{BBCZ05}.  Actually, in this section we prove a more
general evaluation of a finite trigonometric sum for any multiplicative primitive Dirichlet character modulo $m$.

\begin{corollary} \label{cor: multipl character sums}
Choose a positive integer $m$, and let $\tilde{\chi}$ be a primitive Dirichlet character modulo $m$.
Then for any integer $b\in\{1,\ldots,m-1\}$, arbitrary $\alpha,\varphi\in \mathbb{R}$, we have that
\begin{equation} \label{eq: mult char sum cos}
\sum_{j=0}^{m-1}\overline{\tilde{\chi}(j)} \cos^n \left(\frac{2\pi j b}{m} + \alpha\right)=\frac{m}{2^n \tau(\tilde{\chi})}\sum\limits_{r=0}^{m-1}  \sum\limits_{_{\substack{ d \in\{0,\ldots, n\} \\ m\mid (2d-n)b-r }}}\tilde{\chi}(r)  \binom{n}{d} \exp\left(i\alpha(n-2d)\right)
\end{equation}
and
\begin{align} \label{eq: mult char sum sin}
\sum_{j=0}^{m-1}\overline{\tilde{\chi}(j)} \sin^n \left(\frac{2\pi j b}{m} + \varphi \right) = \frac{m}{2^n \tau(\tilde{\chi})} \sum\limits_{r=0}^{m-1}  \sum\limits_{_{\substack{ d \in\{0,\ldots, n\} \\ m\mid (2d-n)b-r }}}\tilde{\chi}(r)  \binom{n}{d} i^{2d-n}\exp\left(i\varphi(n-2d)\right),
\end{align}
where $\tau(\tilde{\chi})$ denotes the Gauss sum associated to the Dirichlet character $\tilde{\chi}$ modulo $m$.
\end{corollary}
\begin{proof}
Let us prove the first identity. The second one follows by taking $\alpha=\varphi -\pi/2$ in  \eqref{eq: mult char sum cos}.

We multiply the identity \eqref{eq: cos twisted sum} (where $\chi(j)=\exp(2 \pi i j/r)$) by $\tilde{\chi}(r)$, and take the sum over all $r\in\{0,\ldots, m-1\}$ to get
\begin{multline*}
\sum_{j=0}^{m-1} \cos^n \left(\frac{2\pi j b}{m} + \alpha \right) \sum_{r=0}^{m-1}\tilde{\chi}(r) e^{2\pi i jr/m} =\frac{m}{2^n}  \sum_{r=0}^{m-1}\tilde{\chi}(r)  \sum\limits_{_{\substack{ d \in\{0,\ldots, n\} \\ m\mid (2d-n)b-r }}}\binom{n}{d} \exp\left(i\alpha(n-2d)\right).
\end{multline*}
Then, a direct application of the formula
$$
\sum_{r=0}^{m-1}\tilde{\chi}(r)e^{2\pi i rj/m}=\overline{\tilde{\chi}(j)}\tau(\tilde{\chi})
$$
completes the proof of the first identity.
\end{proof}

When $\tilde{\chi}$ is an odd, real multiplicative character modulo $m$, then (see Chapter 1 of \cite{BEW98})
$$
 \sum_{r=0}^{m-1}\tilde{\chi}(r) e^{2\pi i jr/m} = i\sqrt{m}\tilde{\chi}(j).
$$
Hence, for odd $n$, equation \eqref{eq: mult char sum sin} becomes the identity
$$
\sum_{j=0}^{m-1}\tilde{\chi}(j) \sin^n \left(\frac{2\pi j b}{m} + \varphi \right) = \frac{\sqrt{m}}{2^n}  \sum_{r=0}^{m-1}\tilde{\chi}(r)  \sum\limits_{_{\substack{ d \in\{0,\ldots, n\} \\ m\mid (2d-n)b-r }}}\binom{n}{d} \exp\left(i\varphi(n-2d)\right) (-1)^{d-(n-1)/2}.
$$
Next, we observe that for an odd $n$ and when $d$ runs through the integers $0,\ldots,n$,
the difference $n-2d$ runs through the odd integers from $-n$ to $n$. Hence, by applying the formula $\exp(ix)+\exp(-ix)=2\cos x$ one easily deduces (under the assumption that $nb<m$) that the right-hand side of the above formula may be expressed as
$$
\frac{\sqrt{m}}{2^{n-1}}  \sum\limits_{_{\substack{ d, r \geq 0 \\ r+2db=nb }}}(-1)^{d-(n-1)/2} \binom{n}{d} \cos\left(\varphi(n-2d)\right) \tilde{\chi}(r),
$$
which proves Theorem 1.2. of \cite{BH10}.
\begin{remark} \rm
By taking $n=1$ and $b=1$ and even/odd real character $\tilde{\chi}$ in Corollary \ref{cor: multipl character sums} one may deduce the Gauss theorems for the evaluation of the Gauss sum $\tau(\tilde{\chi})$.

Namely, when $\tilde{\chi}$ is an even primitive character modulo $m$, Equation \eqref{eq: mult char sum cos} yields the identity
$$
\sum_{j=0}^{m-1}\tilde{\chi}(j) \cos \left(\frac{2\pi j}{m} + \alpha\right) = \frac{m}{\tau(\tilde{\chi})}\cos\alpha.
$$
When $\alpha=0$ the left-hand side of the above equation equals $\tau(\tilde{\chi})$, because we assume $\tilde{\chi}$ is even. Hence, one immediately gets $\tau^2(\tilde{\chi})=m$. If $\tilde{\chi}$ is an odd, primitive character modulo $m$, Equation \eqref{eq: mult char sum cos} becomes
$$
\sum_{j=0}^{m-1}\tilde{\chi}(j) \cos\left(\frac{2\pi j }{m} + \alpha\right) = -\frac{im}{\tau(\tilde{\chi})} \sin\alpha.
$$
When $\alpha=-\pi/2$, using the fact that the character $\tilde{\chi}$ is odd, the left-hand side of the above display reduces to $\frac{1}{i} \tau(\tilde{\chi})$. Therefore, the above identity yields the Gauss formula $\tau^2(\tilde{\chi}) =-m$ for the odd character.

%Applying the Gauss theorem to the sum over even and odd characters we end up with two very elegant identities:
%$$
%\sum_{j=0}^{m-1}\tilde{\chi}(j) \cos \left(\frac{2\pi j}{m} + \alpha\right) = %\sqrt{m}\cos\alpha
%$$
%for even primitive characters modulo $m$, and
%$$
%\sum_{j=0}^{m-1}\tilde{\chi}(j) \cos\left(\frac{2\pi j }{m} + \alpha\right) = %-\sqrt{m}\sin\alpha
%$$
%for odd primitive characters modulo $m$.
\end{remark}

\subsection{Sums which include powers of linear combinations}

Let us now consider an example of a
trigonometric identity for a power of a linear combination of two cosines, which stems from an application of Theorem \ref{thm: main} when $d=2$.

\begin{example} \rm \label{ex: d=2}
Let us take $d=2$, $\mathbf{m}=(m_1,m_2)$ for positive integers $m_1,m_2$ and let $S=\{\pm (1,0),$ $\pm (0,1)\}$ with $\pi_S(s)=1/4$, for all four elements $s\in S$. Let $\bm{\beta}=(\beta_1, \beta_2)\in\mathbb{R}^2$. We apply results of Corollary \ref{cor: main trace f-la} to deduce a trigonometric identity for the sum of powers of a linear combination of cosine functions with shifted arguments.

The eigenvalues $\lambda_{\mathbf{n}}^{\mathbf{m},\bm{\beta}}$ for any $\mathbf{n}=(a,b)$, $a\in\{0,\ldots,m_1-1\}$ and $b\in\{0,\ldots,m_2-1\}$ are given by:
$$
\lambda_{(a,b)}^{(m_1,m_2), (\beta_1,\beta_2)}= \frac{1}{2}\left(\cos\left(2\pi\frac{a+\beta_1}{m_1}\right) + \cos\left(2\pi\frac{b+\beta_2}{m_2}\right)\right).
$$
On the other hand, in view of Proposition \ref{ heat kernel geometric sum}, the right-hand side of \eqref{eq: main formula - cor}, for a fixed positive integer $n$, can be simplified as
$$
4^{-n}m_1m_2\sum\limits_{v=0}^{n/2}\sum_{k_{1}=-\left\lfloor
\frac{2v}{m_1}\right\rfloor }^{\left\lfloor \frac{2v}{m_1}\right\rfloor
}\sum_{k_{2}=-\left\lfloor \frac{(n-2v)}{m_2}\right\rfloor }^{\left\lfloor
\frac{(n-2v)}{m_2}\right\rfloor }e^{2\pi i (\beta_1k_1+\beta_2k_2)}\binom{n}{2v}\binom{2v}{v+k_{1}m_{1}/2}\binom{n-2v}{(n-2v+m_{2}k_{2})/2},
$$
which yields the following trigonometric identity:
\begin{multline*}
\sum_{a=0}^{m_1-1} \sum_{b=0}^{m_2-1} \left(\cos\left(2\pi\frac{a+\beta_1}{m_1}\right) + \cos\left(2\pi\frac{b+\beta_2}{m_2}\right)\right)^n =\\ \frac{m_1m_2}{2^n}\sum\limits_{v=0}^{n/2}\sum_{k_{1}=-\left\lfloor
\frac{2v}{m_1}\right\rfloor }^{\left\lfloor \frac{2v}{m_1}\right\rfloor
}\sum_{k_{2}=-\left\lfloor \frac{(n-2v)}{m_2}\right\rfloor }^{\left\lfloor
\frac{(n-2v)}{m_2}\right\rfloor }e^{2\pi i (\beta_1k_1+\beta_2k_2)}\binom{n}{2v}\binom{2v}{v+k_{1}m_{1}/2}\binom{n-2v}{(n-2v+m_{2}k_{2})/2}.
\end{multline*}
\end{example}

\subsection{Further sums which include powers of linear combinations}\label{further_sums}
Let us take $d=2$ and $S=\{(\pm 1, 0), (0,0), (0,\pm 2)\}$ with probabilities $1/4$ for
the elements $(\pm 1,0)$ and $(0,0)$ and probabilities $1/8$ for the elements $(0,\pm 2)$.  In
this instance, Corollary \ref{cor: main trace f-la} will yield an evaluation of the series
$$
\sum\limits_{k_{1}=0}^{m_{1}-1}
\sum\limits_{k_{2}=0}^{m_{2}-1}\left(\cos\left(\frac{2\pi k_{1}a_{1}}{m_{1}}+ \beta_{1}\right) + \cos^{2}\left(\frac{2\pi k_{2}a_{2}}{m_{2}}+ \beta_{2}\right)\right)^{n}.
$$
In a similar vein, one can consider general $d$ and devise the set $S$ with corresponding probabilities so
that our Corollary \ref{cor: main trace f-la} yields an evaluation of the series
$$
\sum\limits_{k_{1}=0}^{m_{1}-1}\cdots
\sum\limits_{k_{d}=0}^{m_{d}-1}\left(\cos^{h_{1}}\left(\frac{2\pi k_{1}a_{1}}{m_{1}}+ \beta_{1}\right)+\cdots +
\cos^{h_{d}}\left(\frac{2\pi k_{d}a_{d}}{m_{d}}+ \beta_{2}\right)\right)^{n}
$$
for any sequence $\{h_{j}\}$ of positive integers.

By choosing probabilities on $S$ in a suitable way, one can also evaluate weighted multiple trigonometric sums with nonnegative weights. (For an occurrence of such a sum in mathematical modelling, see e.g. \cite{F-DG-D14}.) We will leave the complete development of these formulas
to the interested reader.

Finally, let us note that some multidimensional trigonometric sums can, of course, be written in terms of products of one-dimensional sums and binomial coefficients. However, our approach provides a physical interpretation of this multidimensional sum as a trace of a certain random walk on a graph after $n$ steps; see \cite{OP20} for a related result.

\section{The heat kernel for a product of groups} \label{sec: heat kernel product of groups}

In this section we want to study the heat kernel of the Cayley graph of a
product of finite abelian groups $G_{1}\times G_{2}$ with symmetric subset $%
S_{1}\times S_{2},$ where $S_{i}\subset G_{i}$ are symmetric subsets of
$G_{i}$, and $\pi_i$ are probability distributions in $S_i$, $i=1,2$.  We notice that in this case $S_{1}\times S_{2}$ must also be
symmetric and $\pi_1\times \pi_2$ is a probability distribution in $S_1 \times S_2$.

\subsection{The uniform probability distribution}
We first deal with the case in which each $\pi_i$, $i=1,2$ is the uniform distribution. We will omit it from the notation.
We will denote the Cayley graph $\mathcal{C}(G_{1}\times G_{2},S_{1}\times
S_{2})$ simply as $X_{G_{1}\times G_{2}}.$ Similarly, we will write $%
X_{G_{i}}$ for $\mathcal{C}(G_{i},S_{i})$.
We start by relating the Hilbert space of $L^{2}$ function on $G_{1}\times
G_{2}$ with those of $L^{2}(G_{i},\mathbb{C)}, i=1,2$.

\begin{proposition}	\label{canonical}
There is a canonical isomorphism of Hilbert spaces
	between $L^{2}(G_{1}\times G_{2},\mathbb{C})$ and $L^{2}(G_{1},\mathbb{C)}%
	\otimes L^{2}(G_{2},\mathbb{C)}$ given by%
	\begin{eqnarray*}
		\phi \overset{}{:}L^{2}(G_{1},\mathbb{C)}\otimes _{\mathbb{C}}L^{2}(G_{2},%
		\mathbb{C)} &\rightarrow &L^{2}(G_{1}\times G_{2},\mathbb{C}) \\
		f\otimes g &\mapsto &f\times g.
	\end{eqnarray*}
\end{proposition}

\begin{proof} The proof follows by applying standard arguments, so we omit it here.
%	The natural map
%	\begin{equation*}
%		B:L^{2}(G_{1},\mathbb{C)}\times L^{2}(G_{2},\mathbb{C)}\longrightarrow
%		L^{2}(G_{1}\times G_{2},\mathbb{C})
%	\end{equation*}
%	that sends $(f,g)$ into $f\times g$ is clearly bilinear and thus extends to
%	a unique $\mathbb{C}$-map $\phi $ as above. One readily checks that $\phi $
%	preserves the inner product since
%%%		f,h\right\rangle \left\langle g,l\right\rangle
%	\end{equation*}
%	is equal to $\left\langle f\times g,h\times l\right\rangle .$
%	
%	Let $\delta _{x_{i}}$ be the characteristic functions of $x_{i}\in G_{i}.$%	We know that $\mathfrak{B}_{i}=\{\delta _{x_{i}}:x_{i}\in G_{i}\}$ are bases
%	for $L^{2}(G_{i},\mathbb{C)}$, $i=1,2$. Hence,
%	\begin{equation*}
%		\mathfrak{B}_{1}\mathfrak{\otimes B}_{2}=\{\delta _{x_{1}}\otimes \delta
%		_{x_{2}}:(x_{1},x_{2})\in G_{1}\times G_{2}\}
%	\end{equation*}
%	is a basis for the tensor product space $L^{2}(G_{1},\mathbb{C)}\otimes _{%
%		\mathbb{C}}L^{2}(G_{2},\mathbb{C)}.$ To see that $\phi $ is an isomorphism
%	it is enough to check that $\phi $ sends the basis $\mathfrak{B}_{1}%
%	\mathfrak{\otimes B}_{2}$ injectively onto the set $\{\delta _{x_{1}}\times
%	\delta _{x_{2}}:(x_{1},x_{2})\in G_{1}\times G_{2}\}$, which is a basis for $%
%	L^{2}(G_{1}\times G_{2},\mathbb{C}).$ This is a straightforward
%	computation.
\end{proof}

Let $A_{G_{i}}$ and $A_{G_{1}\times G_{2}}$ be the adjacency operators of $%
X_{G_{i}}$ and $X_{G_{1}\times G_{2}}$, respectively. As one naturally would expect one has $%
A_{G_{1}}\otimes_{\mathbb{C}} A_{G_{2}}\simeq $ $A_{G_{1}\times G_{2}}$ canonically. In
fact, one has the following general result.

\begin{proposition}
	\label{canonical 2}
	There is a canonical isomorphism between $A_{G_{1}\times
		G_{2}}$ and $A_{G_{1}}\otimes A_{G_{2}}$.
\end{proposition}

\begin{proof}
	Consider the diagram%
	\begin{equation*}
		\begin{array}{ccc}
			L^{2}(G_{1},\mathbb{C)}\otimes _{\mathbb{C}}L^{2}(G_{2},\mathbb{C)} &
			\overset{A_{G_{1}}\otimes A_{G_{2}}}{\mathbb{\longrightarrow }} &
			L^{2}(G_{1},\mathbb{C)}\otimes _{\mathbb{C}}L^{2}(G_{2},\mathbb{C)} \\
			\downarrow \phi  &  & \downarrow \phi  \\
			L^{2}(G_{1}\times G_{2},\mathbb{C}) & \overset{A_{G\times G_{2}}}{\mathbb{%
					\longrightarrow }} & L^{2}(G_{1}\times G_{2},\mathbb{C})%
		\end{array}%
	\end{equation*}%
	where the vertical arrows are the isomorphisms defined in Proposition \ref{canonical}.
	We claim this diagram is commutative.
	
	First, the composite map $\phi \circ A_{G_{1}}\otimes A_{G_{2}}$ takes any
	basis vector $\delta _{x_{1}}\otimes \delta _{x_{2}}$ into the function
	given by%
	\begin{equation}
		A_{G1}(\delta _{x_{1}})\times A_{G_{2}}(\delta
		_{x_{2}})(a,b)=\sum\limits_{s_{1}\in S_{1}}\delta
		_{x_{1}}(a+s_{1})\sum\limits_{s_{2}\in S_{2}}\delta _{x_{2}}(b+s_{2})
		\label{je1}
	\end{equation}%
	which is equal to $1,$ if and only if $a+s_{1}=x_{1}$ and $b+s_{2}=x_{2}$.
	
	On the other hand, $A_{G_{1}\times G_{2}}\circ \phi $ takes the basis vector
	$\delta _{x_{1}}\otimes \delta _{x_{2}}$ into the function given by
	\begin{equation}
		A_{G_{1}\times G_{2}}(\delta
		_{(x_{1},x_{2})})(a,b)=\sum\limits_{(s_{1},s_{2})\in S_{1}\times
			S_{2}}\delta _{(x_{1},x_{2})}(a+s_{1},b+s_{2}).  \label{je2}
	\end{equation}%
	From this, it is obvious that (\ref{je1}) and (\ref{je2}) are the same.
	
	This proves the commutativity of the diagram, meaning that
	\begin{equation}
		A_{G_{1}}\otimes A_{G_{2}}=\phi ^{-1}\circ A_{G_{1}\times G_{2}}\circ \phi,
		\label{je3}
	\end{equation}%
	as asserted.
\end{proof}

The previous proposition allows us to relate the heat kernel of $%
X_{G_{1}\times G_{2}}$ with the heat kernels of $G_{i}$, $i=1,2$. These operators
correspond to random walks on $G_{1}\times G_{2}$ and on each $G_{i}$, $i=1,2$. We know that $%
K_{X_{G_{i}}}(x_{i},y_{i};n)=A_{G_i}^{n}(\delta _{x_{i}})(y_{i}).$ From (\ref{je3}%
) we see that%
\begin{equation*}
	(A_{G_{1}}\otimes A_{G_{2}})^{n}=A_{G_{1}}^{n}\otimes A_{G_{2}}^{n}=\phi
	^{-1}\circ A_{G_{1}\times G_{2}}^{n}\circ \phi \, .
\end{equation*}%
Thus,
\begin{equation}
	A_{G_{1}\times G_{2}}^{n}=\phi \circ A_{G_{1}}^{n}\otimes A_{G_{2}}^{n}\circ
	\phi ^{-1}.  \label{je4}
\end{equation}%
On the other hand,
\begin{equation*}
	K_{X_{G_{1}\times G_{2}}}((x_{1},x_{2}),(y_{1},y_{2});n)=A_{G_{1}\times
		G_{2}}^{n}(\delta _{x_{1}}\times \delta _{x_{2}})(y_{1},y_{2}) \,.
\end{equation*}%
By (\ref{je4}) one has
\begin{equation*}
	\left( A_{G_{1}\times G_{2}}^{n}(\delta _{x_{1}}\times \delta
	_{x_{2}})\right) (y_{1},y_{2})=A_{G_{1}}^{n}(\delta _{x_{1}})(y_{1})\times
	A_{G_{2}}^{n}(\delta _{x_{2}})(y_{2}) \,.
\end{equation*}%
Thus,
\begin{equation*}
	K_{X_{G_{1}\times
		G_{2}}}((x_{1},x_{2}),(y_{1},y_{2});n)=K_{X_{G_{1}}}(x_{1},y_{1};n)K_{X_{G_{2}}}(x_{2},y_{2};n) \,.
\end{equation*}

Using a straightforward induction argument, one then obtains the following corollary.

\begin{corollary}
	Let $G_{i}$, $i=1,\ldots ,d$, be finite abelian groups, and let $%
	S_{i}\subset G_{i}$ be symmetric subsets. Then the
	heat kernel of the Cayley graph $\mathcal{C}(G_{1}\times \cdots G_{d},S_{1}\times
	\cdots \times S_{d})$ is the product of the heat kernels of $\mathcal{C}(G_{i},S_{i})$, $i=1,\ldots ,d$.
\end{corollary}

\subsection{Arbitrary probability distribution}

The fact that the heat kernel of a product is the product of the heat
kernels follows from the fact that the adjacency operator of the Cayley
graph $X$ is
canonically identifiable with the adjacency operators of each $%
X_{G_i}.$ This identification readily extends to the general case of
a product of weighted graphs $\mathcal{C}(G_{i},S_{i},\pi _{i})$,  $i=1,\ldots ,d$, where each heat
kernel is taken with a $\chi _{\beta_i}$ -twist.

For this, we first
notice that if $G=G_{1}\times \cdots \times G_{d}$, $S=S_{1}\times \cdots
\times S_{d}$ and $\pi =\pi _{1}\cdots \pi _{d}$ then $X_{G}=\mathcal{C}(G,S,\pi )$ is a
weighted graph since $\pi =\pi _{1}\cdots \pi _{d}$ is a probability
distribution in $S.$

\begin{proposition} \label{prop: product of HK}
	\label{general factors}Let $X_{G_i}=\mathcal{C}(G_{i},S_{i},\pi _{i})$, $i=1,\ldots ,d$ be the
	Cayley graphs of the abelian groups $G_{i}$ with symmetric subset $S_{i}\subset
	G_{i}$ and probability distributions $\pi _{i}$. Let $K_{X_{G_{i}}%
	}(x_i,y_i;n)$ be the discrete time heat kernel of each $X(G_{i},S_{i},\pi _{i}),$
	twisted by $\chi_{\mathbf{\beta}_i}$,  $i=1,\ldots,d$.  Then, the heat kernel on $X_G$, twisted by $\bm{\beta=(\beta_1,\ldots,\beta_d)}$ is given for $x=(x_1,\ldots,x_d)$ and $y=(y_1,\ldots,y_d)$ by
	\begin{equation*}
		K_{X_{G}}(x,y;n)=\prod\limits_{i=1}^{d}K_{X_{G_{i}}}%
		(x_i,y_i;n).
	\end{equation*}
\end{proposition}

\subsection{Multidimensional examples}

We start with the proof of formula \eqref{eq product of d terms}, which is deduced by the following application of Proposition \ref{prop: product of HK}.

\begin{example} \rm

\label{ex: d any}
Let $\mathbf{m}$ denote the $d$-tuple
	$(m_{1},\ldots ,m_{d})$ and  let $G_{i}=\mathbb{Z}_{m_{i}}$ and $G=G_{\mathbf{m}}.$ Let $%
S_{i}=\{\pm 1\}$ so that $S=S_{1}\times \cdots \times S_{d}$ is the set\textrm{\ }%
\begin{equation*}
	S=\{(\epsilon _{1},\ldots ,\epsilon _{d}):\epsilon _{i}=\pm
		1,i=1,\ldots ,d\}.
\end{equation*}%
Let $\pi _{S}(s)=2^{-d}$, for all $s\in S$. Then,
	for any $\mathbf{\beta }\in \mathbb{R}^{d}$, and any $\mathbf{l}=(\ell
	_{1},\ldots ,\ell _{d})\in G_{\mathbf{m}}$ one has
	\begin{equation*}
		\lambda _{\mathbf{l}}^{\mathbf{m},\mathbf{\beta }}=\prod_{k=1}^{d}\cos
		\left( 2\pi (\ell _{k}+\beta _{k})/m_{k}\right) .
	\end{equation*}%
	By Proposition \ref{general factors}, and using Example \ref{ex. d=1 whole group HK}, we see that
	\begin{equation*}
		K_{X_{G}}(\mathbf{k}\mathbf{m},\mathbf{0};n)=2^{-dn}\binom{n}{(n+k_{1}m_{1})/2}\cdots \binom{n}{(n+k_{r}m_{r})/2},
	\end{equation*}%
	under the assumption that $n+k_{j}m_{j}$, $j=1,\ldots ,d$ are even natural
	numbers and $|k_{j}m_{j}|\leq n$, $j=1,\ldots ,d$. In all other cases $K_{X_{G}}(
	\mathbf{k}\mathbf{m},\mathbf{0};n)=0$.

Therefore, application of \eqref{eq: main formula - cor} yields the formula \eqref{eq product of d terms}.
	
\end{example}
The following are two explicit examples of computation of triple sums, using  \emph{Maple}.

\begin{example} \rm
For any $n,m_{1},m_{2},m_{3} \in \mathbb{Z}_{>0}$, let
\begin{equation*}
S(n,m_{1},m_{2},m_{3})=\sum\limits_{a_{1}=0}^{m_{1}-1}\sum\limits_{a_{2}=0}^{m_{2}-1}
\sum\limits_{a_{3}=0}^{m_{3}-1}\cos ^{n}\left(\frac{2\pi a_{1}}{m_{1}}\right)\cos ^{n}\left(\frac{2\pi a_{2}}{m_{2}}\right)
\cos ^{n}\left(\frac{2\pi a_{3}}{m_{3}}\right).
\end{equation*}%
A computation with \emph{Maple} shows that $S(100,40,60,80)=\frac{A}{B},$ where
\begin{eqnarray*}
A &=&3^{6}\times 5^{3}\times 11^{3}\times 13\times 19^{2}\times 29\times
31^{3}\times 83\times 89\times 97^{3} \times 173\times 2699\times 1107114391 \\ &&\times 13231313\times 54570781 \times 60580339\times 20078765421593524568089
\end{eqnarray*}
 and $B=2^{283}.$

For $n=1000$ and $m_1=4,m_2=6,m_3=8$ one gets $S(1000,4,6,8)=\frac{A}{B},$ where
\begin{eqnarray*}
A &=&876866552760850968690689699007021449838100397008270720601894950619774096\\
   &&264083055091427022164624231541328760169885937934982363723536748993854275\\
   &&972950653275500842045793701644136700068694867276069269821412225340930469\\
   &&061982186126624381189674140721945180423650252562631800324801976874415916\\
   &&701971625874609267575978451276158394571564128535848221079129665688994077\\
   &&401903486179240002782019024043632897153590480491978931583944336917869593\\
   &&3770866195399966721
\end{eqnarray*}
and $B=2^{1495}$.  We were unable to compute the prime factorization of $A$.

Let us note that, though it is well known that the sum $S(n,m_{1},m_{2},m_{3})$ is a rational number, numerical computation of the sum becomes quite difficult for large values of $n$ because $n$-th powers of cosines are very small numbers so the precision in the numerical computation must be extremely high.
\end{example}

In our last example we show how to prove \eqref{eq:example 20}.

\begin{example} \rm \label{ex: final d=2}
 Let $X_G$ denote the weighted Cayley graph $G=\mathcal{C}(\mathbb{Z}^2/\mathbf{m}\mathbb{Z}^2,S,\pi_S)$, where $\mathbf{m}=(2m_1,m_2)$, $S=\{(\pm a, \pm b)\}$ is a set of four elements and $\pi_S(s)=1/4$ for all $s\in S$. Then, the heat kernel on $X_G$ twisted by a character $\chi_{\bm{\beta}}$, where $\bm{\beta}=\left(\frac{\alpha_1 m_1}{\pi a}, \frac{m_2(\alpha_2-\pi/2)}{2\pi b} \right)$ equals a product of the heat kernel $K_{X_{G_1}}$ on $G_1=\mathcal{C}(\mathbb{Z}/2m_1\mathbb{Z},S_1,\pi_{S_1})$, where $S_1=\{-a,a\}$, $\pi_{S_1}(s)=1/2$ for all $s\in S_1$ twisted by $\chi_{\frac{\alpha_1 m_1}{\pi a}}$ and the heat kernel $K_{X_{G_2}}$ on $G_2=\mathcal{C}(\mathbb{Z}/m_2\mathbb{Z},S_2,\pi_{S_2})$, where $S_2=\{-b,b\}$, $\pi_{S_2}(s)=1/2$ for all $s\in S_2$ twisted by $\chi_{\frac{m_2(\alpha_2-\pi/2)}{2\pi b}}$.

 Proposition \ref{prop: product of HK} yields that
$$
K_G((m_1,0),(0,0);k)= K_{G_1}(m_1,0;k) K_{G_2}(0,0;k).
$$
Therefore, when combining formulas \eqref{eq: heat k spectral exp}, \eqref{eq: defn of chi,n} and \eqref{eq: defn of lambda,n} in this setting,
we deduce that
\begin{align*}
\sum_{j=0}^{2m_1-1}\sum_{\ell=0}^{m_2-1}(-1)^j&\cos^k\left(\frac{\pi ja}{m_1}+\alpha_1\right)\sin^k\left(\frac{2\pi \ell b}{m_2}+\alpha_2\right)\\&=2m_1m_2\exp(-\pi i \alpha_1 m_1/a)K_{G_1}(m_1,0;k) K_{G_2}(0,0;k).
\end{align*}
Next, we apply the results of Example \ref{ex: d=1, b any beta any} to compute $ K_{G_1}(m_1,0;k)$ and $ K_{G_2}(0,0;k)$ twisted by $\chi_{\frac{\alpha_1 m_1}{\pi a}}$ and  $\chi_{\frac{m_2(\alpha_2-\pi/2)}{2\pi b}}$, respectively. In doing so, we get that
$$
K_{G_1}(m_1,0;k) = 2^{-k}e^{i\alpha_1m_1/a} \sum\limits_{_{\substack{ d_1 \in\{0,\ldots, k\} \\ a(2d_1-k)/m_1-1/2 \in \mathbb{Z} }}}  \binom{k}{d_1} \exp\left(i\alpha_1(k-2d_1)\right),
$$
and
$$
K_{G_2}(0,0;k) = 2^{-k}\sum\limits_{_{\substack{ d_2 \in\{0,\ldots, k\} \\ b(2d_2-k)/m_2 \in \mathbb{Z} }}}  \binom{k}{d_2} \exp\left(i\alpha_2(k-2d_2)\right) i^{2d_2-k}.
$$
With all this, we have proved the explicit evaluation \eqref{eq:example 20}.
%\begin{multline*}
%\sum_{j=0}^{2m_1-1}\sum_{\ell=0}^{m_2-1}(-1)^j\cos^k\left(\frac{\pi ja}{m_1}+\alpha_1\right)\sin^k\left(\frac{2\pi \ell b}{m_2}+\alpha_2\right)= \frac{2m_1m_2}{4^k} \exp\left(ik(\alpha_1+\alpha_2)\right) \cdot\\\cdot \sum\limits_{_{\substack{ d_1,d_2 \in\{0,\ldots, k\} \\ a(2d_1-k)/m_1-1/2 \in \mathbb{Z}, \text{  } b(2d_2-k)/m_2 \in \mathbb{Z} }}} \binom{k}{d_1}  \binom{k}{d_2} \exp\left(-2i(\alpha_1d_1+\alpha_2d_2)\right) i^{2d_2-k}.
%\end{multline*}
\end{example}

\section{Concluding remarks}

%\subsection{Sums of a product of terms with different exponents}

From Proposition \ref{prop: product of HK} and the discussion in Example \ref{ex: d any}
and Example \ref{ex: final d=2} one sees that our methodology applies to yield an explicit
evaluation of the sum of terms of the type
\begin{equation}\label{general_k_sum}
\left(\prod\limits_{j=1}^{d}\cos\left(\frac{2\pi k_{j}a_{j}}{m_{j}}+ \beta_{j}\right)\right)^{k}.
\end{equation}
It is elementary that
%\begin{align*}
%&\left(\frac{\partial^{2}}{\partial \beta_{1}^{2}} + k \right)
%\left(\prod\limits_{j=1}^{d}\cos\left(\frac{2\pi k_{j}a_{j}}{m_{j}}+ %\beta_{j}\right)\right)^{k}
%\\&=k(k-1)\cos^{k-2}\left(\frac{2\pi k_{1}a_{1}}{m_{1}}+ %\beta_{1}\right)\sin^{2}\left(\frac{2\pi k_{1}a_{1}}{m_{1}}+ \beta_{1}\right)
%\left(\prod\limits_{j=2}^{d}\cos\left(\frac{2\pi k_{j}a_{j}}{m_{j}}+ %\beta_{j}\right)\right)^{k}
%\end{align*}
%and that
\begin{equation*}
\left(\frac{\partial^{2}}{\partial \beta_{1}^{2}} + k^{2} \right)
\prod\limits_{j=1}^{d}\cos^k\left(\frac{2\pi k_{j}a_{j}}{m_{j}}+ \beta_{j}\right)
=k(k-1)\cos^{k-2}\left(\frac{2\pi k_{1}a_{1}}{m_{1}}+ \beta_{1}\right)
\prod\limits_{j=2}^{d}\cos^k\left(\frac{2\pi k_{j}a_{j}}{m_{j}}+ \beta_{j}\right).
\end{equation*}
Therefore, from our evaluation of \eqref{general_k_sum}, we can repeatedly apply differential operators
of the form $\partial_{\beta_j}^{2} + c$, for appropriate $c$, and derive an evaluation for sums of
terms of the type
\begin{equation}\label{general_sum}
\prod\limits_{j=1}^{d}\cos^{2h_{j}}\left(\frac{2\pi k_{j}a_{j}}{m_{j}}+ \beta_{j}\right)
\end{equation}
for any set of integers $\{h_{j}\}$.
Furthermore, one can employ the Pythagorean identity in order to obtain sums of terms of the type
\begin{equation} \label{mixed products}
\prod\limits_{j=1}^{d}\cos^{2h_{1,j}}\left(\frac{2\pi k_{j}a_{j}}{m_{j}}+ \beta_{j}\right)
\sin^{2h_{2,j}}\left(\frac{2\pi k_{j}a_{j}}{m_{j}}+ \beta_{j}\right)
\end{equation}
for any set $\{h_{1,j},h_{2,j}\}$ of non-negative integers.  It is not clear how tractable such expressions
may be; nonetheless, the above discussion does provide a guide to a somewhat straightforward algorithm.

Finally, let us describe a means by which one can extend the above assertion to a general
monomial whose exponents are not necessarily even.  To do so, one needs to establish formulas for the sums of terms of the form
\begin{equation}\label{gen_product_one_variable}
\cos^{h_{1}}\left(\frac{2\pi ka}{m}+ \beta\right)\sin^{h_{2}}\left(\frac{2\pi ka}{m}+ \beta\right)
\end{equation}
for general non-negative integers $h_{1}$ and $h_{2}$.  The cases $h_{1}=0$ or $h_{2}=0$ follow
from our main results, as illustrated by Corollary \ref{ex: d=1} with $\beta$ replaced by $m\beta/2\pi$
in the case $h_{2}=0$ and $\beta$ replaced by $m(\beta-\pi/2)/2\pi$ in the case $h_{1}=0$.
Through repeated use of the Pythagorean identity, it suffices to consider the case $h_{2}=1$.
In that instance, such evaluations follow from the case $h_{2}=0$ and that
$$
\frac{\partial}{\partial \beta}
\cos^{h}\left(\frac{2\pi ka}{m}+ \beta\right)
=-h\cos^{h-1}\left(\frac{2\pi ka}{m}+ \beta\right)
\sin\left(\frac{2\pi ka}{m}+ \beta\right).
$$
It remains to be seen if the resulting formulas, though obtainable, are manageable.

In a more direct approach one can take in Proposition \ref{prop: product of HK} a set $S$ equals to a product of $2d$ sets of the form $S_{j}=\{\pm a_j\}$, for $j=1,\ldots,d$ , and $S_{j}=\{\pm b_j\}$, for $j=d+1,\ldots,2d$.  In this case  one obtains a sum of products of terms of the form
\begin{equation} \label{product}
\prod\limits_{j=1}^{d}\cos^{}\left(\frac{2\pi k_{j}a_{j}}{m_{j}} \right)\prod\limits_{j=d+1}^{2d}\cos^{}\left(\frac{2\pi k_{j}b_{j}}{m_{j}} \right)\,.
 \end{equation}
Then one can use  the Chebyshev's polynomials of the first type
\[
T_{n}(x)=\sum\limits_{r=0}^{\left\lfloor n/2\right\rfloor }\frac{(-1)^{r}}{%
n-r}\binom{n-r}{r}(2x)^{n-2r}
\]
to express factors of the form $\cos^n(\alpha)$ in terms of $T_n(\cos(n\alpha))$ and lower powers $\cos^{n-2r}(\alpha)$. Those factors involving sine functions can be obtained by a character twist with $\beta_j=0$, for $j=1,\ldots,d$, and equal to $-m_j/(2\pi)$, for $j=d+1, \ldots,2d$. Then one can proceed in a  recursive manner to compute sums of mixed terms as in (\ref{mixed products}). An explicit algorithm to do this will be developed and implemented in the subsequent article \cite{CHJSV22a}.

%\subsection{Sums of values of polynomials}

%With the discussion from the previous subsection, we arrive at the following conclusion.

%\it
%Let $P(x_{1},\cdots, x_{d};y_{1},\cdots, y_{d})$ be any polynomial in $2d$ variables with complex %coefficients.
%Then, with notation as above, there is a finite step algorithm which explicitly evaluates the series
%\begin{equation}\label{polynomial_sum}
%\sum\limits_{k_{1}=0}^{m_{1}-1} \cdots \sum\limits_{k_{d}=0}^{m_{d}-1}
%P\left(\cdots, \cos\left(\frac{2\pi k_{j}a_{j}}{m_{j}}+ \beta_{j}\right),\cdots;
%\cdots, \sin\left(\frac{2\pi k_{j}a_{j}}{m_{j}}+ \beta_{j}\right),\cdots\right)
%\end{equation}
%as a finite sum involving the degree and coefficients of $P$, binomial coefficients, and exponentials
%which are linear in $\{\beta_{j}\}$.  \rm

%Additional conclusions are possible.  For example, if the coefficients of $P$ are integers
%and one takes $\beta_{j}=0$ for all $j$, then by elementary Galois theory, the value of
%\eqref{polynomial_sum} is a rational number whose denominator is a power of $2$ and whose
%numerator is a multiple of $m_{1}\cdots m_{d}$; see Remark \ref{rem:Galois_argument}.

\vspace{5mm}\noindent
Carlos A. Cadavid \\
Department of Mathematics \\
Universidad Eafit \\
Carrera 49 No 7 Sur-50 \\
Medell\'in, Colombia \\
e-mail: ccadavid@eafit.edu.co

\vspace{5mm}\noindent
Paulina Hoyos  \\
Department of Mathematics \\
The University of Texas at Austin \\
C2515 Speedway, PMA 8.100 \\
Austin, TX 78712
U.S.A. \\
e-mail: paulinah@utexas.edu

\vspace{5mm}\noindent
Jay Jorgenson \\
Department of Mathematics \\
The City College of New York \\
Convent Avenue at 138th Street \\
New York, NY 10031
U.S.A. \\
e-mail: jjorgenson@mindspring.com

\vspace{5mm}\noindent
Lejla Smajlovi\'c \\
Department of Mathematics \\
University of Sarajevo\\
Zmaja od Bosne 35, 71 000 Sarajevo\\
Bosnia and Herzegovina\\
e-mail: lejlas@pmf.unsa.ba

\vspace{5mm}\noindent
Juan D. V\'elez \\
Department of Mathematics \\
Universidad Nacional de Colombia\\
Carrera 65 Nro. 59A - 110\\
Medell\'in, Colombia\\
e-mail: jdvelez@unal.edu.co

\end{document}